\documentclass[12pt,a4paper]{article}

\usepackage{graphicx}
\usepackage{amsfonts}
\usepackage{float}
\usepackage[latin1]{inputenc}
\usepackage{blkarray}
\usepackage{color}
\usepackage{amsmath}
\usepackage{amstext}
\usepackage{amssymb}
\usepackage{array}
\usepackage{multicol}
\usepackage{pdfpages}
\usepackage{multirow}
\usepackage{threeparttable}
\usepackage{url}
\usepackage{enumerate}

\usepackage[ruled,vlined]{algorithm2e}
\usepackage{algorithmic}

\usepackage{subcaption}
\def\N{\Bbb N}
\def\R{\Bbb R}

\def\D{\cal D}

\def\<{\langle}
\def\>{\rangle}

\def\Chi{\raise .3ex \hbox{\large $\chi$}}

\def\n{{\bf }}

\def\[{\Bigl [}
\def\]{\Bigr ]}
\def\({\Bigl (}
\def\){\Bigr )}
\def\[{\Bigl [}
\def\]{\Bigr ]}
\def\({\Bigl (}
\def\){\Bigr )}

\def\div{{\mbox{\rm div}}}
\def\dsp{\displaystyle}

\def\x{{\bf x}}

\def\n{{\bf n}}
\def\s{{\bf s}}

\def\u{{\bf q}}

\def\G{{\Gamma}}

\def\D{{\cal D}}
\def\cells{{\cal M}}
\def\faces{{\cal F}}
\def\nodes{{\cal V}}
\def\edges{{\cal E}}

\def\FG{{{\cal F}_\Gamma}}

\begin{document}

\title{Parallel Vertex Approximate Gradient discretization of hybrid dimensional Darcy flow and transport in discrete fracture networks}
\author{
F. Xing\thanks{Laboratoire de Math\'ematiques J.A. Dieudonn\'e, UMR 7351 CNRS, University Nice Sophia Antipolis, team COFFEE, INRIA Sophia Antipolis M\'editerran\'ee, Parc Valrose 06108 Nice Cedex 02, France, and BRGM Orl\'eans France, feng.xing@unice.fr},         
R. Masson\thanks{Laboratoire de Math\'ematiques J.A. Dieudonn\'e, UMR 7351 CNRS, University Nice Sophia Antipolis, and team COFFEE, INRIA Sophia Antipolis M\'editerran\'ee, Parc Valrose 06108 Nice Cedex 02, France, roland.masson@unice.fr}, 
S. Lopez\thanks{BRGM, scientific and Technical Center, 3 avenue Claude Guillemin, BP 36009, 45060 Orl\'eans Cedex 2 France, s.lopez@brgm.fr}
}
\maketitle


\maketitle

\begin{abstract}
This paper proposes a parallel numerical algorithm to simulate the flow and the transport 
in a discrete fracture network taking into account the mass exchanges with the surrounding matrix. 
The discretization of the Darcy fluxes  is based on the Vertex Approximate Gradient finite volume scheme 
adapted to polyhedral meshes and to heterogeneous anisotropic media, and 
the transport equation is discretized by a first order upwind scheme combined with an Euler explicit 
integration in time. 
The parallelization is based on the SPMD (Single Program, Multiple Data) paradigm and relies 
on a distribution of the mesh on the processes with one layer of ghost cells in order to allow 
for a local assembly of the discrete systems. The linear system for the Darcy flow is solved using different linear solvers and preconditioners implemented in the PETSc and Trilinos libraries. 
The convergence of the scheme is validated on two original analytical solutions with one and four intersecting fractures. 
Then, the parallel efficiency of the algorithm is assessed on up to 512 processes with 
different types of meshes, different matrix fracture 
permeability ratios, and different levels of complexity of the fracture network. 
\end{abstract}

\section{Introduction}

\subsection{Hybrid dimensional flow and transport models}

This article deals with the simulation of the 
Darcy flow and transport in fractured porous media for which the fractures are 
modeled as interfaces of codimension one. In this framework, the $d-1$ dimensional 
flow and transport in the fractures is 
coupled with the $d$ dimensional flow and transport in the matrix leading to the so called   
hybrid dimensional Darcy flow and transport model. 

For the Darcy flow model, we focus on the particular case where the pressure is continuous at 
the interfaces between the fractures and the matrix domain. This type of 
Darcy flow model introduced in \cite{MAE02}, \cite{AKMR05} corresponds physically to pervious 
fractures for which the ratio of the normal permeability of the 
fracture to the width of the fracture is large compared with the ratio of the permeability of 
the matrix to the size of the domain. 
Note that it does not cover the case of fractures acting as barriers 
for which the pressure is discontinuous at the matrix fracture interfaces 
(see \cite{FNFM03}, \cite{MJE05}, \cite{ABH09} for discontinuous pressure models). 
It is also assumed in our model that the pressure is continuous 
at the fracture intersections.
It corresponds  
to the assumption that the ratio between the permeability at the fracture intersections and the width 
of the fracture is large compared to the ratio between the tangential 
permeability of each fracture and its length. 
We refer to \cite{FFSR14} and \cite{SFHW15} for more general reduced models 
taking into account discontinuous pressures at fracture intersections in dimension $d=2$. 

The hybrid dimensional transport model is derived in  \cite{MAE02} in the case of a convection diffusion 
flux for the matrix and fracture concentration. In this work, a purely advective model is considered.
It requires the specification of the transmission conditions at the matrix fracture interfaces and at fracture intersections which, to our knowledge, 
have not been done so far at the continuous level. 

The discretization of the hybrid dimensional Darcy flow model with continuous pressures 
has been the object of several works. In \cite{KDA04} a cell-centred Finite Volume scheme using a 
Two Point Flux Approximation (TPFA) is proposed assuming the orthogonality of the mesh and 
isotropic permeability fields. Cell-centred Finite Volume schemes 
can be extended to general meshes and anisotropic permeability fields 
using  MultiPoint Flux Approximations (MPFA) 
following the ideas of \cite{TFGCH12}, \cite{SBN12},\cite{AELHP152D} and \cite{AELHP153D}. 
In \cite{MAE02} and \cite{HF08} a Mixed Finite Element (MFE) method is proposed, 
and Control Volume Finite Element Methods (CVFE) 
using nodal unknowns have been introduced for such models in 
\cite{BMTA03}, \cite{RJBH06}, \cite{MF07}, \cite{MMB2007}, \cite{G2009}. 
The Hybrid Finite Volume and Mimetic finite difference schemes, belonging  
to the family of Hybrid Mimetic Mixed Methods \cite{DEGH13},  
have been extended to hybrid dimensional models in \cite{FFJR16}, 
\cite{AFSVV16} as well as in \cite{GSDFN}, \cite{BHMS2016} 
in the more general Gradient Discretization framework \cite{DEGGH16}. 
Non-matching discretizations of the fracture and matrix meshes are studied 
in \cite{DS12}, \cite{FS13}, \cite{BPS14} and \cite{SFHW15}. 

Regarding the hybrid dimensional advective transport model, an explicit first order upwind scheme combined 
with the MPFA Darcy fluxes is used in \cite{AELHP152D}, \cite{AELHP153D}, and \cite{SBN12}. 
At fracture intersections, the authors neglect the accumulation term  and the concentration unknown 
is eliminated using the flux conservation 
equation in order to avoid severe restrictions on the time step caused by the small volumes. 
A CVFE method is used in \cite{RJBH06} with a first order upwind approximation and a fully implicit time integration 
of the two phase flow model to avoid small time steps. Higher order methods have also been developed 
in the CVFE method of \cite{MMB2007} using a MUSCL type second order scheme for the saturation equation and also 
in \cite{HF08} where a Discontinuous Galerkin method is used for the transport saturation equation  with an Euler 
implicit time integration in the fracture network and an explicit time integration in the matrix domain. 
In \cite{Haegland2009}, a streamline method is developed in 2D based on the hybrid dimensional Darcy flow velocity field. 
The solution is very accurate for purely advective transport but this method 
requires that the fractures be expanded 
and seems difficult to extend to the case of a complex 3D network in practice.

\subsection{Content and objectives of this work}

In this work, we focus on 
the Vertex Approximate Gradient (VAG) scheme introduced 
in \cite{Eymard.Herbin.ea:2010} for diffusion problems 
and extended in \cite{BGGM14}, \cite{GSDFN},  \cite{BHMS2016}  
to hybrid dimensional Darcy flow models. 
The VAG scheme uses nodal and fracture-face unknowns 
in addition to the cell unknowns which can be eliminated without any fill-in. 
Thanks to its essentially nodal feature, it leads to a sparse discretization on tetrahedral or mainly tetrahedral meshes. 
The VAG scheme has the major advantage, compared 
with the CVFE methods  of \cite{BMTA03}, \cite{RJBH06}, \cite{MF07} or \cite{MMB2007}, 
that it avoids the mixing of the control volumes at the fracture matrix interfaces, 
which is a key feature for its coupling with the transport model. As shown 
in \cite{BGGM14} for two phase flow problems, the VAG scheme allows 
for a coarser mesh size at the matrix fracture interface for a given accuracy. 
For the discretization of the transport hybrid dimensional model, 
we will use in this work a simple first order upwind scheme with explicit time integration. 
The extension to second order MUSCL type discretization will be considered in a future work. 
Our main objective in this paper is to develop a parallel algorithm 
for the VAG discretization of hybrid dimensional Darcy flow and transport models, 
and to assess the parallel scalability of the algorithm. 

Starting from the hybrid dimensional Darcy flow model of \cite{BGGM14} and \cite{GSDFN}, 
we first derive the hybrid dimensional transport model for a general fracture network 
taking into account fracture intersections and the coupling with the matrix domain. 
Then, the VAG discretization of the Darcy flow model is recalled and 
the VAG Darcy fluxes are used to discretize the transport model with an upwind first order discretization in space 
and an Euler explicit time integration. A key feature of this discretization is the definition of the control volumes which 
is adapted to the heterogeneities of the porous medium. This can be achieved thanks to the fact that,   
on the one hand, the VAG scheme keeps the cell unknowns and, on the other hand, the VAG Darcy fluxes are 
constructed independently of the definition of the control volumes. In particular, the control volumes are constructed 
in such a way that, at matrix fracture interfaces, the volume is taken only in the fracture. Otherwise, the fracture 
will be enlarged artificially and the front velocity will not be accurately approximated in the fractures 
as it it the case for usual CVFE methods. 
Note also that we do not eliminate the concentration 
unknowns at fracture intersections as was done in \cite{SBN12}, \cite{AELHP152D} and \cite{AELHP153D} for cell centred discretizations. 
In the case of a nodal discretization like the VAG scheme, this elimination is not possible since 
these unknowns are connected to the matrix and it is not 
needed since the size of the control volumes at fracture intersections 
is roughly the same as 
the size of any control volume located at the matrix fracture interface. 

Our parallelization of the hybrid dimensional flow and transport numerical model 
is based on the SPMD (Single Program, Multiple Data) paradigm. It 
relies on a distribution of the mesh on the processes with one layer of ghost cells in order to allow 
for a local assembly of the discrete systems. 
The linear system for the Darcy flow is solved using different linear solvers 
and preconditioners implemented in the PETSc and Trilinos libraries. 

In order to validate the convergence of the scheme, two analytical solutions 
are constructed for the hybrid dimensional 
flow and transport model. 
We consider the case of a single non-immersed fracture as well as the case of 
four intersecting fractures. The analytical solutions for the transport model are  obtained by integration of the matrix and 
fracture equations along the characteristics of the velocity field taking into account source terms coming from the 
matrix fracture transmission conditions. 
Then, we study the parallel scalability of the Darcy flow and transport solvers on up to 512 processes. 
Our numerical investigation includes different levels of complexity of the fracture network with a number of fractures ranging from a few to a few hundreds. 
It covers different types of meshes namely hexahedral, tetrahedral and prismatic meshes 
as well as a large range of permeability ratios between the fracture network and the matrix domain. 
In addition, the influence of the choices of the linear solver and of the preconditioner 
is also studied for the solution of the Darcy flow equation. \\

The paper is organized as follows. Section \ref{sec_model} recalls the geometrical and functional 
setting introduced in \cite{GSDFN} for a general 2D fracture network immersed in a surrounding 3D matrix. 
Then, the hybrid dimensional Darcy flow and transport models are introduced. 
In Section \ref{sec_VAG}, the VAG discretization is recalled for the Darcy flow model and extended to the transport model. The parallel implementation of the scheme is detailed in section \ref{sec_parallel}. 
Section \ref{sec_tests} is devoted to the description of the test cases including the analytical solutions and to 
the numerical investigation of the parallel scalability of the algorithm.

\section{Hybrid dimensional Darcy Flow and Transport Model in Fractured Porous Media}
\label{sec_model}
\subsection{Discrete Fracture Network and functional setting}

Let $\Omega$ denote a bounded domain of $\R^d$, $d=2,3$ 
assumed to be polyhedral for $d=3$ and polygonal for $d=2$. 
To fix ideas the dimension will be fixed to $d=3$ when it needs to be specified, 
for instance in the naming of the geometrical objects or for the space discretization 
in the next section. The adaptations to the case $d=2$ are straightforward. 

We consider the asymptotic model introduced in \cite{MAE02} 
where fractures are represented as interfaces of codimension 1. 
Let $I$ be a finite set and let 
$
\overline \Gamma = \bigcup_{i\in I} \overline \Gamma_i
$  
and its interior $\Gamma = \overline \Gamma\setminus \partial\overline\Gamma$ 
denote the network of fractures $\Gamma_i\subset \Omega$, $i\in I$, such that each $\Gamma_i$ is 
a planar, polygonal, simply connected, open domain included in an oriented plane ${\cal P}_i$ of $\R^d$. 
It is assumed that the angles of $\Gamma_i$ 
are strictly smaller than $2\pi$ and that $\Gamma_i\cap\overline\Gamma_j=\emptyset$ for all $i\neq j$. 
For all $i\in I$, 
let us set $\Sigma_i = \partial\Gamma_i$, $\Sigma_{i,j}= \Sigma_i\cap\Sigma_j$, $j\in I\setminus\{i\}$, 
$\Sigma_{i,0} = \Sigma_i\cap\partial\Omega$,
$\Sigma_{i,N} = \Sigma_i\setminus(\bigcup_{j\in I\setminus\{i\}}\Sigma_{i,j}\cup \Sigma_{i,0})$, 
and $\Sigma = \bigcup_{(i,j)\in I\times I, i\neq j} \Sigma_{i,j}$. 
It is assumed that $\Sigma_{i,0} = \overline\Gamma_i\cap\partial\Omega$. 
\begin{figure}[h!]
\begin{center}
\includegraphics[width=0.5\textwidth]{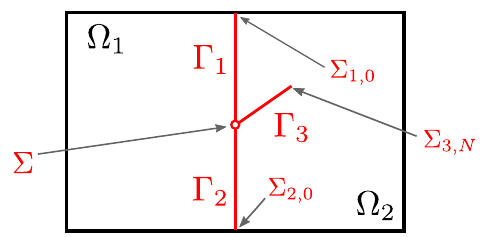}   
\caption{Example of a 2D domain with 3 intersecting fractures $\Gamma_1, \Gamma_2, \Gamma_3$ 
and 2 connected components $\Omega_1$, $\Omega_2$.}
\label{fig_network}
\end{center}
\end{figure}
We will denote by $d\tau({\bf x})$ the $d-1$ dimensional Lebesgue measure on $\Gamma$.
On the fracture network $\G$, we define the function space
$
L^2(\Gamma) = \{v = (v_i)_{i\in I}, v_i\in L^2(\Gamma_i), i\in I\}, 
$
endowed with the norm $\|v\|^2_{L^2(\Gamma)} = \sum_{i\in I} \|v_i\|^2_{L^2(\Gamma_i)}$. 
Its subspace $H^1(\Gamma)$ is defined as the space of functions $v = (v_i)_{i\in I}$ such that $v_i\in H^1(\Gamma_i)$, $i\in I$ with continuous traces at the fracture intersections. 
The space $H^1(\Gamma)$ is 
endowed with the norm $\|v\|^2_{H^1(\Gamma)} = \sum_{i\in I} \|v_i\|^2_{H^1(\Gamma_i)}$ and its subspace with vanishing traces on $\Sigma_0 = \bigcup_{i\in I} \Sigma_{i,0}$ is denoted by $H^1_{\Sigma_0}(\Gamma)$. 

Let us also consider the trace operator 
$\gamma_i$ from $H^1(\Omega)$ to $L^2(\Gamma_i)$ as well 
as the trace operator $\gamma$ from $H^1(\Omega)$ to $L^2(\Gamma)$ such that 
$(\gamma v)_i = \gamma_i(v)$ for all $i\in I$. 
 
On $\Omega$, the gradient operator 
from $H^1(\Omega)$ to $L^2(\Omega)^d$ is denoted by $\nabla$. On the fracture 
network $\G$, the tangential gradient $\nabla_\tau$ acting from 
$H^1(\Gamma)$ to $L^2(\Gamma)^{d-1}$ is defined by 
$$
\nabla_\tau v = (\nabla_{\tau_i} v_i)_{i\in I},  
$$
where, for each $i\in I$, the tangential gradient $\nabla_{\tau_i}$ 
is defined from $H^1(\G_i)$ to $L^2(\G_i)^{d-1}$ by fixing a 
reference Cartesian coordinate system of the plane ${\cal P}_i$ containing $\G_i$. 
We also denote by $\div_{\tau_i}$ the divergence operator from 
$H_{\div}(\Gamma_i)$ to $L^2(\Gamma_i)$. \\ 

The function spaces arising in the variational formulation 
of the hybrid dimensional Darcy flow model are
$$
V = \{v\in H^1(\Omega) \mbox{ such that } \gamma v \in H^1(\Gamma)\}, 
$$ 
and its subspace
$$
V^0 = \{v  \in H^1_0(\Omega) \mbox{ such that } \gamma v \in H^1_{\Sigma_0}(\Gamma) \}. 
$$
The space $V^0$ is endowed with the following Hilbertian norm 
$$
\|v\|_{V^0} = \(\|\nabla v\|^2_{L^2(\Omega)^d} 
+ \|\nabla_\tau \gamma v\|^2_{L^2(\Gamma)^{d-1}} \)^{1/2}. 
$$
Let $\Omega_\alpha, \alpha\in {\cal A}$ denote the connected components of 
$\Omega\setminus\overline\Gamma$, with ${\cal A}$ being the set of connected components of $\Omega\setminus\overline\Gamma$. 
Let us define the space $H_{\div}(\Omega\setminus\overline\Gamma)=
\{{\bf q}_m = ({\bf q}_{m,\alpha})_{\alpha\in {\cal A}} \,|\, {\bf q}_{m,\alpha}\in H_{\div}(\Omega_\alpha)\}$. 
Using the orientation of ${\cal P}_i$
we can define the two sides $\pm$ of the fracture $\Gamma_i$, for all $i\in I$.
For all ${\bf q}_m\in H_{\div}(\Omega\setminus\overline\Gamma)$, let  
$\gamma_{\n,i}^\pm {\bf q}_{m}$ denote the normal trace of ${\bf q}_m$  on the side $\pm$ of $\Gamma_i$ with 
the normal oriented outward from the side $\pm$.  
Let us define the Hilbert function space 
$$
\left.\begin{array}{r@{\,\,}c@{\,\,}l}
H(\Omega,\Gamma) =  \displaystyle 
\{&&
{\bf q}_m = ({\bf q}_{m,\alpha})_{\alpha\in {\cal A}}, \, {\bf q}_f = ({\bf q}_{f,i})_{i\in I} \,|\, 
{\bf q}_{m}\in H_{\div}(\Omega\setminus\overline\Gamma), \\
&&\displaystyle {\bf q}_{f,i}\in L^2(\Gamma_i)^{d-1}, \, 
\div_{\tau_i}({\bf q}_{f,i})  - \gamma_{\n,i}^+ {\bf q}_m - \gamma_{\n,i}^- {\bf q}_m \in L^2(\Gamma_i), i\in I
\}, 
\end{array}\right.
$$
and its closed Hilbert subspace 
\begin{eqnarray}
\label{def_W}
\left.\begin{array}{r@{\,\,}c@{\,\,}l}
&&W(\Omega,\Gamma) =  \displaystyle \{
({\bf q}_m, {\bf q}_f)\in  H(\Omega,\Gamma) \,|\, 
\displaystyle \sum_{\alpha\in {\cal A}}\int_{\Omega_\alpha} ({\bf q}_{m,\alpha} \cdot \nabla v 
+ \div({\bf q}_{m,\alpha}) v)d{\bf x} \\
&&\displaystyle + \sum_{i\in I}\int_{\Gamma_i} ({\bf q}_{f,i} \cdot \nabla_{\tau_i}\gamma_i v  
+  (\div_{\tau_i}({\bf q}_{f,i})  - \gamma_{\n,i}^+ {\bf q}_m - \gamma_{\n,i}^- {\bf q}_m) \gamma_i v)d\tau(\x) = 0 
\,\forall\, v \in V^0
\}.  
\end{array}\right.
\end{eqnarray}
The last definition corresponds to imposing in a weak sense the conditions 
$\sum_{i\in I} \gamma_{\n,\Sigma_i} {\bf q}_{f,i}  =  0$ on  $\Sigma\setminus\Sigma_0$  and  
$ \gamma_{\n,\Sigma_i} {\bf q}_{f,i} = 0$ on $\Sigma_{i,N},i\in I$, 
where $\gamma_{\n,\Sigma_i}$ is the normal trace operator on $\Sigma_i$ (tangent to $\Gamma_i$) 
with the normal oriented outward from $\Gamma_i$,  and using the extension 
of $\gamma_{\n,\Sigma_i} {\bf q}_{f,i}$ by zero on $\Sigma\setminus\Sigma_i$.

\subsection{Hybrid dimensional Darcy Flow Model}

In the matrix domain $\Omega\setminus\overline\G$ (resp. in the fracture network $\G$), let us denote 
by $\Lambda_m\in L^{\infty}(\Omega)^{d\times d}$ (resp. $\Lambda_f\in L^{\infty}(\Gamma)^{(d-1)\times (d-1)}$) 
the permeability tensor such that 
there exist  $\overline\lambda_m\geq \underline\lambda_m > 0$ (resp. $\overline\lambda_f\geq \underline\lambda_f > 0$) 
with 
$$
\underline\lambda_m|\boldsymbol{\xi}|^2 \leq (\Lambda_m(\x)\boldsymbol{\xi},\boldsymbol{\xi}) \leq \overline\lambda_m|\boldsymbol{\xi}|^2 
\mbox{ for all } \boldsymbol{\xi} \in \R^d, \x\in \Omega,
$$
(resp. $
\underline\lambda_f|\boldsymbol{\xi}|^2 \leq (\Lambda_f(\x)\boldsymbol{\xi},\boldsymbol{\xi}) \leq \overline\lambda_f|\boldsymbol{\xi}|^2 
\mbox{ for all } \boldsymbol{\xi} \in \R^{d-1}, \x\in\Gamma$). 

We also denote by $\mu$ the fluid viscosity which is assumed constant and by $d_f \in L^\infty(\G)$ 
the width of the fractures assumed to be such that there exist 
${\overline d}_f\geq {\underline d}_f > 0$ with 
$
{\underline d}_f \leq d_f(\x) \leq {\overline d}_f
$
for all $\x\in\Gamma$.

Given $\bar u\in V$, the strong formulation of the hybrid dimensional Darcy flow model amounts to: find $u\in V$ 
and $(\u_m, \u_f) \in W(\Omega,\Gamma)$ such that $u-\bar u\in V^0$ and 
\begin{eqnarray}
\label{modeleCont}
\left\{\begin{array}{r@{\,\,}c@{\,\,}ll}
\div(\u_{m,\alpha}) &=& 0 &\mbox{ on } \Omega_\alpha, \alpha\in {\cal A},\\
\u_{m,\alpha} &=& -{\Lambda_m \over \mu}\nabla u &\mbox{ on } \Omega_\alpha, \alpha\in{\cal A},\\  
\div_{\tau_i}(\u_{f,i})  - \gamma_{\n,i}^+ {\bf q}_m - \gamma_{\n,i}^- {\bf q}_m  &=& 0 &\mbox{ on } \G_i, i\in I, \\
\u_{f,i} &=& -d_f ~{\Lambda_f \over \mu} \nabla_{\tau_i} \gamma_i u  &\mbox{ on } \G_i, i\in I.    
\end{array}\right.
\end{eqnarray}

The weak formulation of \eqref{modeleCont} amounts to: 
find $u\in V$  such that $u-\bar u\in V^0$ and the following 
variational equation is satisfied for all $v\in V^0$: 
\begin{eqnarray}
\label{formVar}
\dsp \int_\Omega  {\Lambda_m(\x)\over \mu}\nabla u(\x) \cdot\nabla v(\x) d\x  
+ \dsp\int_\G  d_f(\x) {\Lambda_f(\x)\over \mu}  \nabla_\tau \gamma u(\x) \cdot\nabla_\tau \gamma v(\x) d\tau(\x)  = 0. 
\end{eqnarray}
The existence and uniqueness of the solution to \eqref{formVar} derives from the Lax Milgram theorem and 
a Poincar\'e inequality stated in \cite{GSDFN}. 

\subsection{Hybrid dimensional transport model}

Let $\gamma_\n$ be the normal trace operator on $\partial\Omega$ with the normal oriented outward from $\Omega$. 
Let us define $\partial\Omega^- = \{\x \in \partial\Omega \,|\, \gamma_\n {\bf q}_{m} (\x) < 0\}$,   
$\Sigma_{i,0}^- = \{\x \in  \Sigma_{i,0} \,|\, \gamma_{\n,\Sigma_i} {\bf q}_{f,i}(\x) < 0\}$, $i\in I$, as well as 
the following subset of $\Sigma\setminus\Sigma_0$:  
$$
\Sigma^- = \{\x \in \Sigma\setminus\Sigma_0 \,|\, \sum_{i\in I}| \gamma_{\n,\Sigma_i}{\bf q}_{f,i}(\x)|  \neq 0\}. 
$$
We consider a linear, purely advective model with velocity $\u_m$ in the matrix domain and 
$\u_f$ in the fracture network. The matrix concentration is denoted by $c_m$ ($c_{m,\alpha}$ in each connected component 
$\Omega_\alpha$, $\alpha\in {\cal A}$) and the fracture concentration, 
representing the average concentration in the fracture width, 
is denoted by $c_f$ ($c_{f,i}$ in each fracture $\Gamma_i$, $i\in I$). 
The 2D equation in the fracture network 
is as usual obtained by integration of the 3D advection equation in the width of the fractures. 
For a purely advective equation, the transmission condition  
at the matrix fracture interfaces states that the input normal flux in the matrix is obtained using the upwind 
fracture concentration $c_f$. 
At the fracture intersection $\Sigma^-$, an additional unknown $c_{f,\Sigma^-}$ must be introduced 
and the transmission conditions state that the normal fluxes sum to zero and that the input normal 
fluxes are obtained using the upwind concentration $c_{f,\Sigma^-}$. \\

Let be given the input boundary conditions 
$\bar c_{m} \in L^\infty(\partial\Omega^-)$, 
$\bar c_{f,i}\in  L^\infty(\Sigma_{i,0}^-)$, $i \in I$, and the initial conditions 
$c_{m}^0 \in L^{\infty}(\Omega\setminus\overline\Gamma)$, $c_{f}^0 \in L^{\infty}(\G)$. 
Let us denote by $\phi_m(\x)$ the porosity in the matrix and by $\phi_f(\x)$ the porosity in the 
fracture network. The transport hybrid dimensional model amounts to:  
find $c_{m} \in L^{\infty}\((\Omega\setminus\overline\Gamma)\times (0,T)\)$, 
$c_{f}\in  L^{\infty}\(\G\times (0,T)\)$, 
and $c_{f,\Sigma^-}\in L^{\infty}\(\Sigma^-\times (0,T)\)$, 
such that: 
\begin{eqnarray}
\label{modeleTracer}
\left\{\begin{array}{r@{\,\,}c@{\,\,}ll}
\phi_m \partial_t c_{m,\alpha} + \div(c_{m,\alpha} {\bf q}_{m,\alpha}) &=& 0 &\mbox{ on } 
\Omega_\alpha\times (0,T), \alpha\in {\cal A}\\
\phi_f d_f  \partial_t c_{f,i} + \div_{\tau_i}(c_{f,i} {\bf q}_{f,i}) &=& 
\gamma_{\n,i}^+ c_m \u_m + \gamma_{\n,i}^- c_m \u_m &\mbox{ on } \G_i\times (0,T), i\in I,\\
(\gamma_{\n,i}^\pm c_m \u_m)^- &=& c_f (\gamma_{\n,i}^\pm \u_m)^-   &\mbox{ on } \G_i\times (0,T), i\in I,\\
(\gamma_{\n,\Sigma_i} c_{f,i} \u_{f,i})^- &=& c_{f,\Sigma^-} (\gamma_{\n,\Sigma_i} \u_{f,i})^- & 
\mbox{ on } (\Sigma_i\setminus\Sigma_{i,0}) \times (0,T), i\in I,\\
\dsp \sum_{j\in I} \gamma_{\n,\Sigma_j} c_{f,j} {\bf q}_{f,j}  &=&  0 & \mbox{ on } (\Sigma\setminus\Sigma_0)\times (0,T),\\
(\gamma_\n c_{m} \u_m)^- &=& \bar c_{m} (\gamma_\n q_m)^-  & \mbox{ on } \partial\Omega \times (0,T), \\
(\gamma_{\n,\Sigma_i} c_{f,i} \u_{f,i})^- &=&\bar c_{f,i}(\gamma_{\n,\Sigma_i} \u_{f,i})^-  & \mbox{ on } \Sigma_{i,0}\times (0,T), i\in I,\\
c_{m} &=& c_{m}^0 & \mbox{ on } (\Omega\setminus\overline\Gamma) \times \{t=0\},\\
c_{f} &=& c_{f}^0 & \mbox{ on } \G\times \{t=0\},
\end{array}\right.
\end{eqnarray}
where the notations $a^+ =\max(a,0)$ and $a^- = \min(a,0)$ are used for all $a\in \R$.

\section{Vertex Approximate Gradient Discretization (VAG)}
\label{sec_VAG}

\subsection{VAG discretization of the Darcy flow model}

In the spirit of \cite{Eymard.Herbin.ea:2010}, we consider generalized polyhedral meshes of $\Omega$ in the sense that the edges of the mesh are linear but its faces are not necessarily planar. Roughly speaking, each face 
is assumed to be defined by the union of the triangles joining each edge of the face to a  
so called face centre. This definition has the advantage to include 
in particular hexahedral cells with non planar faces. 

Let $\cells$ be the set of cells which are disjoint open polyhedral subsets of $\Omega$ such that
$\bigcup_{K\in\cells} \overline{K} = \overline\Omega$. For each $K\in\cells$, 
it is assumed that there exists 
${\x}_K\in K\setminus\partial K$, the so-called ``centre'' of the cell $K$, such that $K$ is star-shaped with respect to ${\x}_K$. We then denote by $\faces_K$ the set of interfaces of non zero $d-1$ dimensional measure among the interior faces
$\overline{K}\cap\overline{L}$, $L\in\cells\setminus\{K\}$, and the boundary interface $\overline{K}\cap\partial\Omega$, which possibly splits in several boundary faces. Let us denote by
$$
\faces = \bigcup_{K\in \cells}\faces_K
$$
the set of all faces of the mesh. Note that the faces are not assumed to be planar, hence the term ``generalized polyhedral mesh''. For $\sigma\in \faces$, let $\edges_\sigma$ be the set of interfaces of non zero, $d-2$ dimensional measure among the interfaces
$\overline{\sigma}\cap\overline{\sigma}'$, $\sigma'\in{\faces}\setminus\{\sigma\}$. Then, we denote by 
$$
\edges = \bigcup_{ \sigma \in \faces }\edges_\sigma
$$
the set of all edges of the mesh. 
Let $\nodes_\sigma = \bigcup_{\left(e,e'\right) \in \edges^2_\sigma, e\neq e'}\big( e\cap e'\big)$ be the set of nodes of $\sigma$. 
For each $K\in \cells$ we define 
$
\nodes_K = \bigcup_{\sigma\in \faces_K} \nodes_\sigma,
$ 
and we also denote by 
$$
\nodes = \bigcup_{K\in \cells} \nodes_K
$$ 
the set of all nodes of the mesh. It is then assumed that for each face $\sigma\in\faces$, there exists a so-called ``centre'' of the face ${{\bf x}}_\sigma \in {\sigma}\setminus \bigcup_{e\in \edges_\sigma} e$ such that
$
{\x}_\sigma = \sum_{\s\in \nodes_\sigma} \beta_{\sigma,\s}~\x_\s, \mbox{ with }
\sum_{\s\in \nodes_\sigma} \beta_{\sigma,\s}=1,
$
and $\beta_{\sigma,\s}\geq 0$ for all $\s\in \nodes_\sigma$; moreover
the face $\sigma$ is assumed to be defined  by the union of the triangles
$T_{\sigma,e}$ defined by the face centre ${\x}_\sigma$
and each edge $e\in\edges_\sigma$. 

The mesh is also supposed to be conforming w.r.t. the fracture network $\G$ in the sense that for each  $i\in I$ 
there exists a subset $\faces_{\G_i}$ of $\faces$ such that 
$
\overline \G_i = \bigcup_{\sigma\in\faces_{\G_i}} \overline{\sigma}.  
$
We will denote by $\faces_\G$ the set of fracture faces $\bigcup_{i\in I} \faces_{\G_i}$. 
The following notations will be used for convenience: 
$$
\cells_\s = \{K\in \cells\,|\, \s\in \nodes_K\}, 
$$
$$
\cells_{\sigma} = \{K\in \cells\,|\, \sigma\in \faces_K\}, 
$$
$$
\faces_{\G,\s} = \{\sigma \in \faces_\Gamma \,|\, \s\in \nodes_\sigma\}, 
$$
and 
$$
\faces_{\Gamma,K} = \faces_K\cap \FG. 
$$
This geometrical discretization of $\Omega$ and $\G$ is denoted in the following by $\D$. \\

The VAG discretization was introduced in \cite{Eymard.Herbin.ea:2010} for diffusive problems 
on heterogeneous, anisotropic media. Its extension to the hybrid dimensional 
Darcy model is based on the following vector space of unknowns: 
$$
X_\D =\{v_K, v_\s, v_\sigma\in\R, K\in \cells, \s\in \nodes, \sigma\in \faces_\G\},
$$
and its subspace with homogeneous Dirichlet boundary conditions on $\partial\Omega$:
$$
X_\D^0 =\{v\in X_\D 
\,|\, v_\s = 0 \mbox{ for } \s\in \nodes_{ext}\}.
$$
where $\nodes_{ext}=\nodes\cap \partial\Omega$ denotes the set of boundary nodes, 
and $\nodes_{int}=\nodes\setminus \partial\Omega$ denotes the set of interior nodes. \\ 

\begin{figure}[H]
\centering
\includegraphics[width=0.3\textwidth]{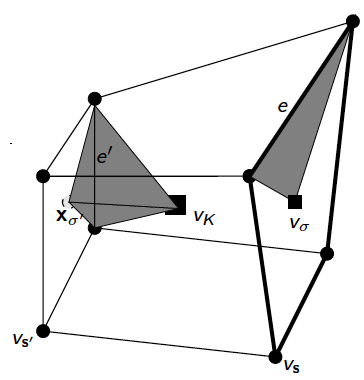}\\
\caption{For a cell $K$ with one fracture face $\sigma$ in bold: 
cell unknown  $v_K$ located at $\x_K$, fracture-face unknown $v_\sigma$ located at $\x_\sigma$, 
node unknowns $v_\s, v_{\s'}$, face centre $\x_{\sigma'}$ of face $\sigma'$, 
triangle $T_{\sigma,e}$ (convex hull of $e$ and $\x_\sigma$), 
triangle $T_{\sigma',e'}$ (convex hull of $e'$ and $\x_{\sigma'}$) 
and tetrahedron $T_{K,\sigma',e'}$ (convex hull of $\x_K$, $\x_{\sigma'}$ and $e'$). }
\label{figddl}
\end{figure}

A finite element discretization of $V$ is built using a tetrahedral sub-mesh of $\cells$ and 
a second order interpolation at the face centres $\x_\sigma$,  
$\sigma\in\faces\setminus\faces_\G$ defined by 
the operator $I_\sigma: X_\D\rightarrow \R$ such that
$$
I_\sigma (v) = \sum_{\s\in \nodes_\sigma} \beta_{\sigma,\s} v_\s. 
$$
The tetrahedral sub-mesh is defined by 
${\cal T} = \{T_{K,\sigma,e}, e\in \edges_\sigma, \sigma\in \faces_K, K\in\cells\}$ 
where  $T_{K,\sigma,e}$ is the tetrahedron joining the cell centre
${\x}_K$ to the triangle $T_{\sigma,e}$ (see Figure \ref{figddl} for examples of 
such tetrahedra).

For a given $v\in X_\D$, we define the function $\pi_{\cal T}v\in V$ as the continuous
piecewise affine function on each
tetrahedron of ${\cal T}$ 
such that $\pi_{\cal T}v({\x}_K) = v_K$, 
$\pi_{\cal T}v(\s) = v_\s$, 
$\pi_{\cal T}v({\x}_\sigma) = v_\sigma$, 
and $\pi_{\cal T}v({\x}_{\sigma'}) = I_{\sigma'}(v)$ 
for all $K\in \cells$, 
$\s\in \nodes$, $\sigma\in \faces_\G$, and $\sigma'\in\faces\setminus\faces_\G$. 
The nodal basis of this finite element discretization will be denoted by
$\eta_K$, $\eta_\s$, $\eta_\sigma$, for $K\in \cells$, $\s\in \nodes$, $\sigma\in \faces_\G$. \\

The VAG discretization of the hybrid dimensional Darcy flow model \eqref{modeleCont} is based on its weak formulation \eqref{formVar}. 
Given $\bar u_\s, \s\in \nodes_{ext}$, it amounts to: find $u_\D\in X_\D$ with $u_\s = \bar u_\s$ for all 
$\s\in {\cal V}_{ext}$ and such that 
for all $v_\D\in X_\D^0$ one has 
\begin{eqnarray}
\label{varformVAG}
\int_\Omega {\Lambda_m(\x)\over \mu} \nabla \pi_{\cal T}u_\D (\x)\cdot \nabla \pi_{\cal T}v_\D(\x) d\x + 
\int_\Gamma d_f(\x) {\Lambda_f(\x)\over \mu} \nabla_\tau \gamma \pi_{\cal T}u_\D(\x)
\cdot\nabla_\tau \gamma \pi_{\cal T}v_\D(\x) d\tau(\x) = 0.  
\end{eqnarray} 

Following \cite{GSDFN}, this Galerkin Finite Element formulation \eqref{varformVAG} 
can be reformulated in terms of discrete conservation laws using the following 
definition of the VAG fluxes. 
For all $v_\D\in X_\D$, the VAG matrix fluxes connect the cell $K\in\cells$
to its nodes and fracture faces $ \nu\in \nodes_K \cup \faces_{\Gamma,K}$: 
\begin{eqnarray}
\label{VAGFluxMatrix}
F_{K,\nu}(v_\D) =  -\int_K {\Lambda_m(\x) \over \mu} \nabla \pi_{\cal T}v_\D(\x)\cdot\nabla\eta_\nu(\x) d\x = \sum_{\nu'\in \nodes_K \cup \faces_{\Gamma,K}}
a_{K,\nu}^{\nu'} (v_K - v_{\nu'}) 
\end{eqnarray} 
with 
$
a_{K,\nu}^{\nu'} = \int_K {\Lambda_m({\x})\over \mu} \nabla \eta_\nu({\x}) \cdot  \nabla\eta_{\nu'}({\x}) d\x. 
$
The VAG fracture fluxes connect the face $\sigma\in\FG$ to its nodes $\s\in\nodes_\sigma$: 
\begin{eqnarray}
\label{VAGFluxFracture}
F_{\sigma,\s}(v_\D) =  -\int_\sigma d_f(\x) {\Lambda_f(\x)\over \mu} \nabla_\tau \gamma \pi_{\cal T}v_\D(\x)
\cdot\nabla_\tau \gamma \eta_\s(\x) d\tau(\x) = \sum_{\s'\in \nodes_\sigma}
a_{\sigma,\s}^{\s'} (v_\sigma - v_{\s'}) 
\end{eqnarray} 
with 
$
a_{\sigma,\s}^{\s'} = \int_\sigma d_f(\x) {\Lambda_f ({\x})\over \mu} \nabla_\tau \gamma \eta_\s({\x}) \cdot  \nabla_\tau \gamma \eta_{\s'}({\x}) 
d\tau(\x).  
$ \\
\begin{figure}[h!]
\begin{center}
\includegraphics[width=0.4\textwidth]{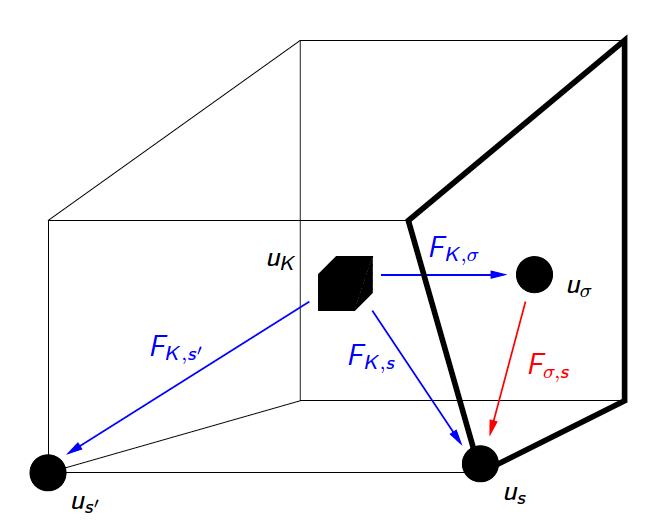} 
\caption{Matrix fluxes (in blue) and fracture fluxes (in red) 
inside a cell $K$ with a 
fracture face $\sigma$ (in bold). The matrix fluxes $F_{K,\nu}$ 
connect the cell $K$ to its nodes and fracture faces $\nu\in \nodes_K \cup \faces_{\Gamma,K}$. 
The fracture fluxes $F_{\sigma,\s}$ connect 
the face $\sigma$ to the nodes $\s\in \nodes_{\sigma}$ of $\sigma$.}
\label{fig_fluxes}
\end{center}
\end{figure}

Then, the Galerkin Finite Element formulation \eqref{varformVAG} is equivalent to: find $u_\D\in X_\D$ 
satisfying the following set of 
discrete conservation equations and Dirichlet boundary conditions: 
$$
\left\{\begin{array}{r@{\,\,}c@{\,\,}ll}
&& \dsp \sum_{\s\in \nodes_K} F_{K,\s}(u_\D) + \sum_{\sigma\in \faces_{\Gamma,K}} F_{K,\sigma}(u_\D) = 0
,  \,\, K\in \cells\\
&&  \dsp \sum_{\s\in {\nodes}_\sigma} F_{\sigma,\s}(u_\D) + \sum_{K\in \cells_\sigma} -  F_{K,\sigma}(u_\D) = 0, 
\,\, \sigma \in {\faces}_\Gamma\\
&&\dsp \sum_{K\in \cells_\s} -  F_{K,\s}(u_\D)
+ \sum_{\sigma \in \faces_{\Gamma,\s}} - F_{\sigma,\s}(u_\D) = 0, \,\, \s \in \nodes_{int},\\
&& u_\s = \bar u_\s, \,\, \s\in  \nodes_{ext}. 
\end{array}\right.
$$

\subsection{First order upwind discretization of the transport model}

\subsubsection{Definition of control volumes}

\begin{figure}[H]
\begin{center}
\includegraphics[width=0.4\textwidth]{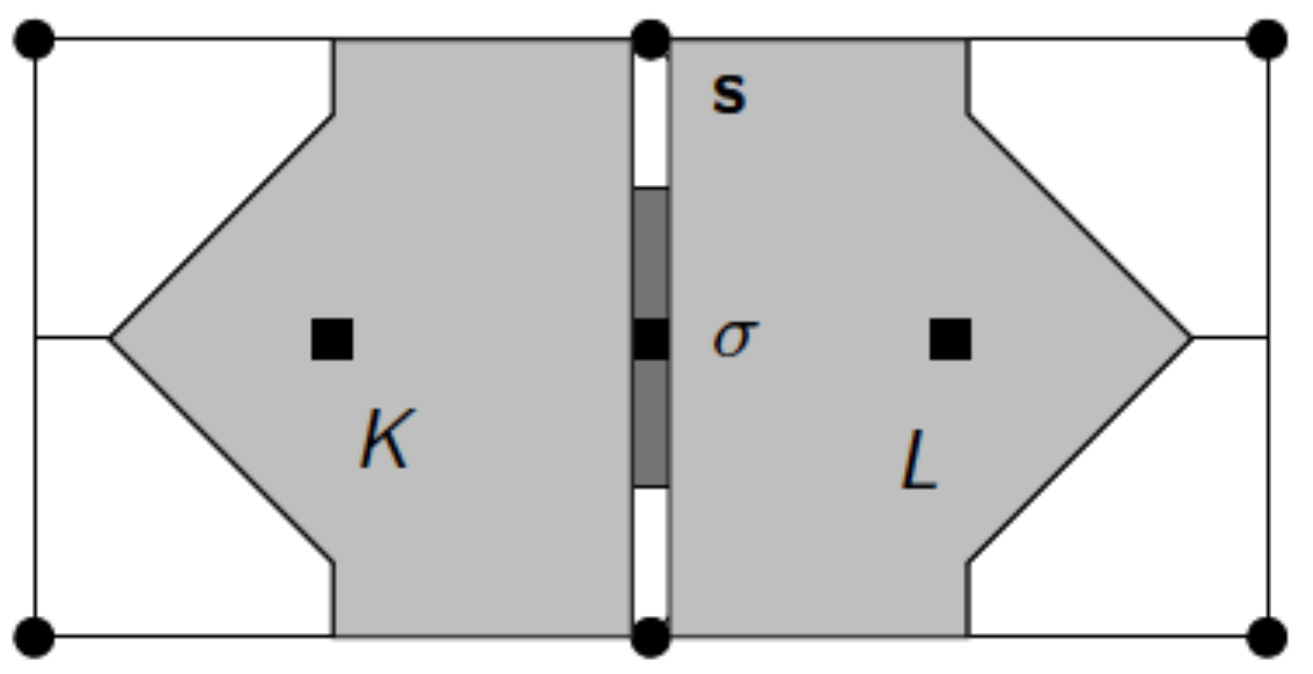}
\includegraphics[width=0.4\textwidth]{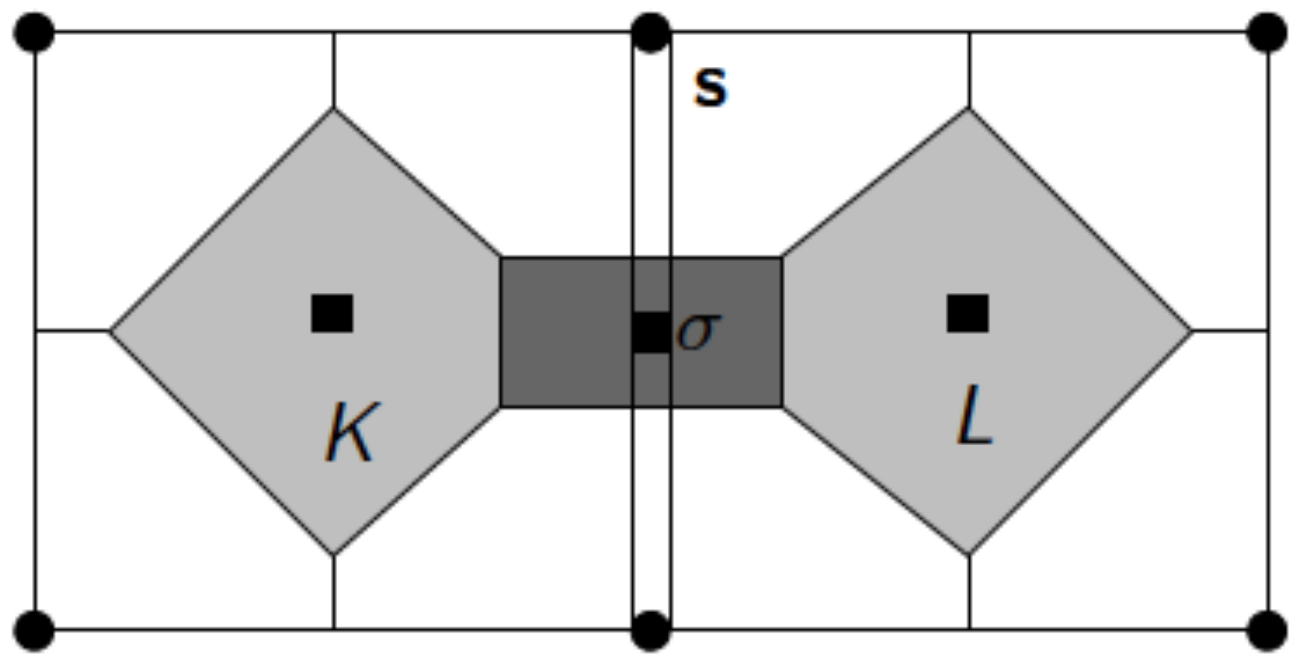}
\caption{Example of choices of the control volumes at cells, fracture face, and nodes, 
in the case of two cells $K$ and $L$ splitted by 
one fracture face $\sigma$ (the width of the fracture has been enlarged in this figure).
The left figure exhibits the good choice with no mixing of fracture and matrix rocktypes while the right figure 
exhibits the bad choice enlarging artificially the fracture.}
\label{fig_vag12}
\end{center}
\end{figure}

The VAG discretization of the hybrid dimensional transport model combines the VAG matrix 
 and fracture  fluxes  \eqref{VAGFluxMatrix}, \eqref{VAGFluxFracture} with the following definition of 
the control volumes based on partitions of the cells and of the fracture faces. 
These partitions are respectively denoted, for all $K\in\cells$,  by 
$$
\overline{K} ~ =  ~\overline{\omega}_K ~ \cup ~ \( \bigcup_{\s\in\nodes_K\setminus(\nodes_{ext}\cup\nodes_\Gamma)} \overline{\omega}_{K,\s} \)
$$
and, for all $\sigma\in\faces_\G$, by 
$$
\overline{\sigma}~=~ \overline{\omega}_\sigma~\cup ~ 
\( \bigcup_{\s\in\nodes_\sigma\setminus\nodes_{ext}} \overline{\omega}_{\sigma,\s} \). 
$$
Then, the control volumes are defined by $\omega_K$ for all cells $K\in \cells$, by $\omega_\sigma$ for all 
fracture faces $\sigma\in \faces_\Gamma$, and  by 
$$
\overline \omega_\s = \bigcup_{K \in \cells_\s} \overline{\omega}_{K,\s}, 
$$
for all nodes $\s\in \nodes_{int} \setminus \nodes_\Gamma$, and by 
$$
\overline \omega_\s = \bigcup_{\sigma \in \faces_{\Gamma,\s}} \overline{\omega}_{\sigma,\s}, 
$$
for all nodes $\s\in \nodes_\Gamma\setminus \nodes_{ext}$. 
Note that this definition avoid the mixing of the fracture and matrix rocktypes at the control volumes 
$\s\in \nodes_\Gamma \setminus \nodes_{ext}$ and $\sigma\in \faces_\Gamma$. This is exhibited in 
Figure \ref{fig_vag12} in comparison with an alternative 
choice mixing the matrix and fracture rocktypes which artificially enlarges 
the fractures. We refer to \cite{GSDFN} for numerical comparisons on a two phase flow model 
of these two types of choices of the control volumes. 
The same idea is applied for all nodes located at different rocktype interfaces. 
In practice, for such a node $\s\in \nodes_{int}\setminus \nodes_\Gamma$ (resp. $\s\in \nodes_\Gamma\setminus\nodes_{ext}$), 
the set $\omega_{K,\s}$ (resp. $\omega_{\sigma,\s}$) should be non empty only for 
the cell(s) $K$ (resp. fracture face(s) $\sigma$) with the largest permeability among those around the node $\s$  
(see  \cite{EHGM-CG-12} for details). \\

In practice, the above partitions of the cells and fracture faces does not need to be built. 
It is sufficient to define the matrix volume fractions 
$$
\alpha_{K,\s} = {\int_{\omega_{K,\s}} d\x \over \int_K d\x}, \s\in\nodes_K\setminus(\nodes_{ext}\cup\nodes_\Gamma), K\in\cells, 
$$
constrained to satisfy $\alpha_{K,\s}\geq 0$, 
and $\sum_{\s\in\nodes_K\setminus(\nodes_{ext}\cup\nodes_\Gamma)}\alpha_{K,\s} \leq 1$, as well as the fracture volume fractions 
$$
\alpha_{\sigma,\s} = {\int_{\omega_{\sigma,\s}} d_f(\x) d\tau(\x) \over \int_\sigma d_f(\x) d\tau(\x)}, 
\s\in\nodes_\sigma\setminus\nodes_{ext}, \sigma\in\faces_\G, 
$$
such that $\alpha_{\sigma,\s}\geq 0$, and $\sum_{\s\in\nodes_\sigma\setminus\nodes_{ext}}\alpha_{\sigma,\s}\leq 1$.  
Then, the porous volumes of the control volumes are set to 
$$
\phi_K = (1 -  \sum_{\s\in \nodes_K \setminus (\nodes_\G\cup\nodes_{ext})}\alpha_{K,s}) \int_K \phi_m(\x)d\x,\,\, K\in\cells
$$
$$
\phi_\sigma = (1 - \sum_{\s \in {\nodes}_\sigma \setminus \nodes_{ext}} \alpha_{\sigma,\s} ) d_{f,\sigma}
\int_\sigma \phi_f(\x)  d\tau(\x),\,\, \sigma\in \faces_\G,
$$
$$
\phi_\s =    \sum_{\sigma\in \faces_{\Gamma,\s}}\alpha_{\sigma,\s} d_{f,\sigma}\int_\sigma \phi_f(\x)   d\tau(\x),\,\, \s \in \nodes_\G\setminus\nodes_{ext}, 
$$
$$
\phi_\s =    \sum_{K\in \cells_{\s}} \alpha_{K,\s} \int_K \phi_m(\x)d\x,\,\, \s \in \nodes \setminus( \nodes_{ext} \cup \nodes_\G), 
$$
with $d_{f,\sigma} = \dsp {\int_\sigma d_f(\x) d\tau(\x) \over \int_\sigma d\tau(\x)}$. 
\subsubsection{Time integration}

For $N\in\N^*$, let us consider the time discretization 
$t^0= 0 < t^1 <\cdots < t^{n-1} < t^n \cdots < t^N = T$ 
of the time interval $[0,T]$. We denote the time steps by 
$\Delta t^n = t^{n+1}-t^{n}$ for all $n=0,\cdots,N-1$. \\

Given  $\bar c_\s, \s\in \nodes_{ext}$ with arbitrary values on  the set of ouput boundary nodes 
$$
\nodes_{ext}^+ = \{\s\in \nodes_{ext} \,|\, F_{K,\s}(u_\D) \geq 0 \, \forall K\in \cells_\s \mbox{ and }  
F_{\sigma,\s}(u_\D) \geq 0 \, \forall 
\sigma \in \faces_{\G,\s} \}, 
$$
and $c_\D^0 \in X_\D$ such that $c_\s^0 = \bar c_\s$ for all $\s\in \nodes_{ext}$, the transport discrete model amounts to 
find $c_\D^{n+1} \in X_\D$ for all $n=0,\cdots,N-1$ satisfying the following discrete conservation laws and Dirichlet input conditions 
$$
\left\{\begin{array}{r@{\,\,}c@{\,\,}ll}
&& \dsp \phi_K {c_{K}^{n+1}- c_{K}^{n} \over \Delta t^n}  
+ \sum_{\s\in {\nodes}_K} H_{K,\s}(c_\D^n) + \sum_{\sigma\in {\faces}_{\Gamma,K}} H_{K,\sigma}(c_\D^n)   = 0,  
\,\, K \in \cells,\\
&&  \dsp  \phi_\sigma {c_{\sigma}^{n+1}- c_{\sigma}^{n} \over \Delta t^n} + 
\sum_{\s \in {\nodes}_\sigma} H_{\sigma,\s}(c_\D^n) 
- \sum_{K \in \cells_\sigma} H_{K,\sigma}(c_\D^n) = 0,
\,\, \sigma \in {\faces}_\Gamma,\\
&&\dsp \phi_\s {c_{\s}^{n+1}- c_{\s}^{n} \over \Delta t^n} -  \sum_{K\in \cells_\s} H_{K,\s}(c_\D^n)
- \sum_{\sigma \in \faces_{\Gamma,\s}} H_{\sigma,\s}(c_\D^n) = 0,
\,\, \s \in \nodes_{int},\\
&& c_\s^{n+1} = \bar c_\s, \,\, \s \in  \nodes_{ext},  
\end{array}\right.
$$
with the following explicit upwind two point fluxes
\begin{equation}
\begin{array}{l}
\label{darcyfluxH}
H_{K,\nu}(c^n_\D) = c_{K}^{n}F_{K,\nu}(u_\D)^+ + c_{\nu}^{n}F_{K,\nu}(u_\D)^- \\ \\
H_{\sigma,\s}(c^n_\D) = c_{\sigma}^n F_{\sigma,\s}(u_\D)^+ + c_\s^n F_{\sigma,\s}(u_\D)^-. 
\end{array}
\end{equation}
The solution of this explicit upwind scheme classically satisfies the following maximum principle 
$$
m \leq c_\mu^{n+1} \leq M \mbox{ for all } \, \mu\in \nodes\cup\faces_\Gamma\cup\cells\setminus {\cal V}_{ext}^+,
$$
with 
$$
M = \max_{\mu\in \nodes\cup\faces_\Gamma\cup\cells\setminus {\cal V}_{ext}^+} c_\mu^0 \, \mbox{ and } \, 
m = \min_{\mu\in \nodes\cup\faces_\Gamma\cup\cells\setminus {\cal V}_{ext}^+} c_\mu^0, 
$$ 
provided that the following Courant-Friedrichs-Lewy (CFL) condition 
\begin{equation}
\label{eq_CFL}
\Delta t^n \leq \min(\Delta t_\cells, \Delta t_{\faces_\Gamma}, \Delta t_\nodes), 
\end{equation}
is satisfied with 
$$
\left\{\begin{array}{r@{\,\,}c@{\,\,}l}
\dsp \Delta t_\cells &=& \dsp \min_{K\in \cells} {\phi_K \over \sum_{\s\in \nodes_K} F_{K,\s}(u_\D)^+ +
\sum_{\sigma\in \faces_{\Gamma,K}} F_{K,\sigma}(u_\D)^+}, \\
\dsp \Delta t_{\faces_\Gamma} &=& \dsp \min_{\sigma\in \faces_\Gamma} 
{\phi_\sigma \over \sum_{\s\in \nodes_\sigma} F_{\sigma,\s}(u_\D)^+  + \sum_{K\in \cells_\sigma} (-F_{K,\sigma}(u_\D))^+}, \\
\dsp \Delta t_{\nodes} &=& \dsp \min_{\s\in \nodes_{int}} 
{\phi_\s \over \sum_{K\in \cells_\s} (-F_{K,\s}(u_\D))^+  + \sum_{\sigma\in \faces_{\Gamma,\s}} (-F_{\sigma,\s}(u_\D))^+}. 
\end{array}\right.
$$

\section{Parallel implementation in ComPASS}
\label{sec_parallel}

The hybrid dimensional Darcy flow and transport discrete model is implemented 
in the framework of the code ComPASS (Computing Parallel Architecture to Speed up Simulations) \cite{ComPass2013}, 
which focuses on parallel high performance simulation (distributed memory, Message Parsing Interface - MPI) adapted to general unstructured polyhedral meshes 
(see \cite{ComPassTwoPhase}). 

\subsection{Mesh non overlapping and overlapping decompositions}

Let us denote by $N_p$ the number of MPI processes. The set of cells $\cells$ is partitioned into $N_p$ subsets $\cells^p, p=1,...,N_p$ using the library METIS \cite{citemetis}. The partitioning of the set of nodes $\mathcal{V}$ and of the set of fracture faces $\mathcal{F}_{\Gamma}$ is defined as follows:
assuming we have defined a global index of the cells $K\in \cells$ 
let us denote 
by $K(\s), \s\in\nodes$ (resp. $K(\sigma)$, $\sigma\in \faces_\Gamma$) the cell with the smallest global index among those of $\cells_\s$ 
(resp. $\cells_\sigma$). Then we set 
$$
\nodes^p = \{\s\in\nodes \,|\, K(\s)\in \cells^p\}, 
$$
and 
$$
\faces_\Gamma^p = \{\sigma\in\faces_\Gamma \,|\, K(\sigma)\in \cells^p\}. 
$$

\begin{figure}[H]
\centering
\includegraphics[width=0.4\textwidth]{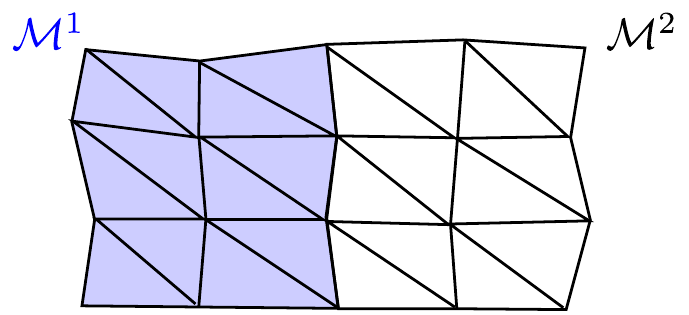}
\hspace{0.05\textwidth}
\includegraphics[width=0.355\textwidth]{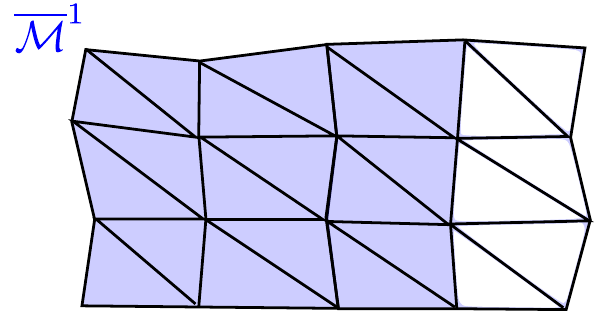} \\
\hspace{0.05\textwidth}
\includegraphics[width=0.4\textwidth]{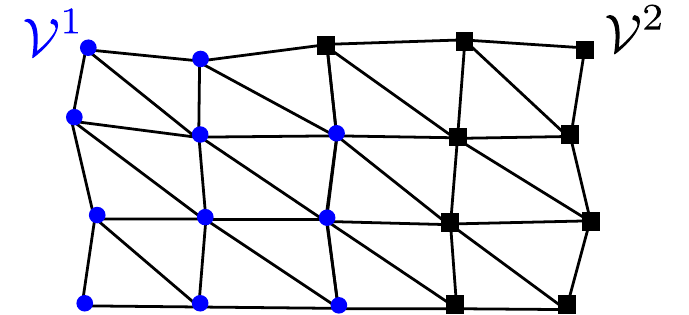}
\hspace{0.04\textwidth}
\includegraphics[width=0.41\textwidth]{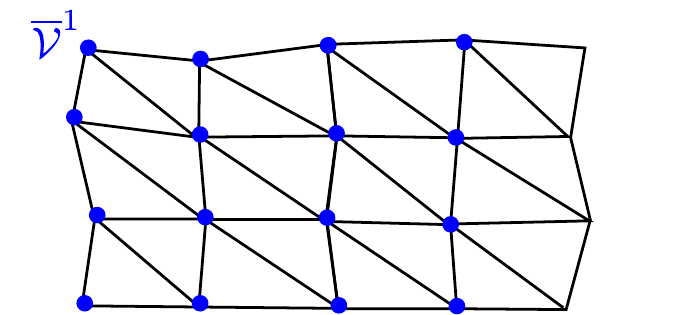} 
\caption{Example of mesh decomposition.}
\label{fig_cellghost} 
\end{figure}

The overlapping decomposition of $\cells$ into the sets
$$
\overline{\cells}^p, \,\, p=1,...,N_p, 
$$
is chosen in such a way that any compact finite volume scheme such as the VAG scheme can be assembled locally on each process. 
Hence, as exhibited in Figure \ref{fig_cellghost}, $\overline{\cells}^p$ is defined as the set of cells sharing a node with a cell of $\cells^p$. 
The overlapping decompositions of the set of nodes and of 
the set of fracture faces follow from this definition:
$$
\overline\nodes^p = \bigcup_{K\in \overline\cells^p} \nodes_K, \,\, p=1,\cdots,N_p, 
$$
and 
$$
\overline\faces_\Gamma^p = \bigcup_{K\in \overline\cells^p} \faces_K\cap\faces_\Gamma, \,\, p=1,\cdots,N_p.  
$$

The partitioning of the mesh is performed by the master process (process 1), and then, each local mesh is distributed to its process.
Therefore, each MPI process contains the local mesh 
($\overline{\cells}^p$, $\overline\nodes^p$,  $\overline\faces_\Gamma^p$), $p=1,2,...,N_p$ which
is splitted into two parts: 
\begin{equation*}
\begin{split}
& \text{own mesh: }  (\cells^p, \nodes^p, \faces_\Gamma^p), \\
& \text{ghost mesh: }  (\overline{\cells}^p \backslash \cells^p, \overline\nodes^p \backslash \nodes^p, \overline\faces_\Gamma^p \backslash \faces_\Gamma^p).
\end{split}
\end{equation*}
We now turn to the parallel implementation of the discrete hybrid dimensional Darcy flow model \eqref{modeleCont} 
and transport model \eqref{modeleTracer}. 

\subsection{Parallelization of the discrete hybrid dimensional Darcy flow}
On each process $p=1,...,N_p$, the local problem of the discrete hybrid dimensional Darcy flow \eqref{modeleCont} is defined 
by the set of unknowns $u_\mu$, $\mu\in \overline\nodes^p\cup \overline\faces_\Gamma^p\cup\overline \cells^p$ and 
the set of equations 
\begin{equation}
\label{equowndarcy}
\left\{\begin{array}{r@{\,\,}c@{\,\,}ll}
&& \dsp \sum_{\s\in {\nodes}_K} F_{K,\s}(u_\D) + \sum_{\sigma\in {\faces}_{\Gamma,K}} F_{K,\sigma}(u_\D) = 0
,  \,\, K\in \overline{\cells}^p, \\
&&  \dsp \sum_{\s\in {\nodes}_\sigma} F_{\sigma,\s}(u_\D) + \sum_{K\in \cells_\sigma} -  F_{K,\sigma}(u_\D) = 0, 
\,\, \sigma \in \mathcal{F}_\Gamma^p, \\
&&\dsp \sum_{K\in \cells_s} -  F_{K,\s}(u_\D)
+ \sum_{\sigma \in \faces_{\Gamma,\s}} - F_{\sigma,\s}(u_\D) = 0, \,\, \s \in \nodes_{int} \cap \mathcal{V}^p,\\
&& u_\s = \bar u_\s, \,\, \s\in  \nodes_{ext} \cap \mathcal{V}^p. 
\end{array}\right.
\end{equation}
Note that this includes the equations of the nodes $\s \in \mathcal{V}^p$, of the fracture faces $\sigma \in \mathcal{F}_\Gamma^p$ 
and of the cells $ K\in \overline{\cells}^p$, both those own cells in $K\in {\cells}^p$ and 
the ghost cells $K\in \overline{\cells}^p\setminus {\cells}^p$. 
The set of equations can be rewritten as the following rectangular linear system 
\begin{equation}
\label{sysowndarcy}
\begin{pmatrix}
A_{v v}^p & A_{v f}^p& A_{v c}^p \\
A_{f v}^p & A_{ff}^p& A_{fc}^p \\
A_{cv}^p & A_{cf}^p& A_{cc}^p
\end{pmatrix}
\begin{pmatrix}
\overline U_v^p \\
\overline U_f^p \\
\overline U_c^p
\end{pmatrix}
=
\begin{pmatrix}
b_v^p \\
b_f^p \\
\overline b_c^p
\end{pmatrix}
\end{equation} 
where $\overline U_v^p \in \mathbb{R}^{\overline{\mathcal{V}}^p} $, $\overline U_f^p \in \mathbb{R}^{\overline{\mathcal{F}}_\Gamma^p}$ 
and $\overline U_c^p \in \mathbb{R}^{\overline{\cells}^p}$ denote the vector of process $p$ own and ghost unknowns at nodes, fracture faces and cells respectively.  
The above matrices have the following sizes 
\begin{eqnarray*}
& A_{vv}^p    \in \mathbb{R}^{ \mathcal{V}^p \times \overline{\mathcal{V}}^p }, \ 
A_{vf}^p \in \mathbb{R}^{ \mathcal{V}^p \times \overline{\mathcal{F}}_{\Gamma}^p }, \
A_{vc}^p      \in \mathbb{R}^{ \mathcal{V}^p \times \overline{\cells}^p }, \\ 
& A_{fv}^p    \in \mathbb{R}^{ \mathcal{F}_{\Gamma}^p \times \overline{\mathcal{V}}^p }, \ 
A_{ff}^p \in \mathbb{R}^{ \mathcal{F}_{\Gamma}^p \times \overline{\mathcal{F}}_{\Gamma}^p }, \
A_{fc}^p      \in \mathbb{R}^{ \mathcal{F}_{\Gamma}^p \times \overline{\cells}^p }, \\ 
& A_{cv}^p    \in \mathbb{R}^{ \overline{\cells}^p \times \overline{\mathcal{V}}^p }, \ 
A_{cf}^p \in \mathbb{R}^{ \overline{\cells}^p \times \overline{\mathcal{F}}_{\Gamma}^p }, \
A_{cc}^p      \in \mathbb{R}^{ \overline{\cells}^p \times \overline{\cells}^p }.
\end{eqnarray*}
and $b_v^p \in \mathbb{R}^{{\mathcal{V}}^p}$, $b_f^p \in \mathbb{R}^{{\mathcal{F}}_\Gamma^p}$, $\overline b_c^p \in \mathbb{R}^{\overline{\cells}^p}$ 
denote the corresponding right hand side vectors. The matrix $A_{cc}^p$ is a non singular diagonal matrix and the cell unknowns can be 
easily eliminated without fill-in leading to the following Schur complement system 
\begin{equation}
\label{eqschur}
{A}^p  
\begin{pmatrix}
\overline U_v^p \\
\overline U_f^p
\end{pmatrix}
=
{b}^p
\end{equation}
with 
\begin{equation*}
A^p:= 
\begin{pmatrix}
A_{vv}^p & A_{vf}^p \\
A_{fv}^p & A_{ff}^p
\end{pmatrix}
-
\begin{pmatrix}
A_{vc}^p \\
A_{fc}^p
\end{pmatrix}
(A_{cc}^p)^{-1}
\begin{pmatrix}
A_{cv}^p & A_{cf}^p
\end{pmatrix},
\
b^p:=
\begin{pmatrix}
b_v^p \\
b_f^p
\end{pmatrix}
-
\begin{pmatrix}
A_{vc}^p \\
A_{fc}^p
\end{pmatrix}
(A_{cc}^p)^{-1} \overline b_c^p,
\end{equation*}
and
\begin{equation}
\label{eqschurcell}
\overline U_c^p = (A_{cc}^p)^{-1} (b_c^p - A_{cv}^p \overline U_v^p - A_{cf}^p \overline U_f^p).
\end{equation}

The linear system \eqref{eqschur} is built locally on each process $p$ and
transfered to the parallel linear solver library PETSc \cite{citepetsc} or Trilinos \cite{citetrilinos}.
The parallel matrix and the parallel vector in PETSc or Trilinos are stored in a distributed manner, i.e. each process stores its own rows. We construct the following parallel linear system 
\begin{equation}
\label{parallellin}
A U= b,
\end{equation}
with
\begin{equation*}
A:=
\begin{pmatrix}
A^1 R^1 \\
A^2 R^2 \\
\vdots \\
A^{N_p} R^{N_p}
\end{pmatrix}
\begin{array}{l}
\big\} \text{ process 1}  \\
\big\} \text{ process 2}  \\
\quad \quad \vdots \\
\big\} \text{ process } N_p
\end{array}, 
\quad
U:=
\begin{pmatrix}
U_v^1 \\
U_f^1 \\
U_v^2 \\
U_f^2 \\
\vdots
\end{pmatrix}
\begin{array}{l}
\bigg\} \text{ process 1} \\
\bigg\} \text{ process 2} \\
\quad \quad \vdots
\end{array}, \quad
b:=
\begin{pmatrix}
b^1 \\
b^2 \\
\vdots \\
b^{N^p}
\end{pmatrix}
\begin{array}{l}
\big\} \text{ process 1 }  \\
\big\} \text{ process 2 }  \\
\quad \quad \vdots \\
\big\} \text{ process } N_p
\end{array}
\end{equation*}
where $R^p,p=1,2,...,N_p$ is a restriction matrix satisfying 
$$R^p U = \begin{pmatrix}
\overline U_v^p \\
\overline U_f^p
\end{pmatrix}.$$
The matrix $A^p R^p$, the vector $\begin{pmatrix}
U_v^p \\
U_f^p 
\end{pmatrix}$
and the vector $b^p$ are stored in process $p$.

The parallel linear system \eqref{parallellin} is solved using 
the GMRES or BiCGStab algorithm preconditioned by different type of preconditioners as discussed in the numerical section. 
The solution of the linear system provides on each process $p$ the solution vector $\begin{pmatrix}
U_v^p \\
U_f^p 
\end{pmatrix}$
of own node and fracture-face unknowns. 
Then, the ghost node unknowns $u_\mu$, $\mu\in (\overline\nodes^p \backslash \nodes^p)$ 
and the ghost fracture-face unknowns $u_\mu$, $\mu\in  (\overline\faces_\Gamma^p \backslash \faces_\Gamma^p)$ are recovered 
by a synchronization step with MPI communications. This synchronization is efficiently 
implemented using a PETSc or Trilinos matrix vector product 
\begin{equation}
\label{commS}
\overline U = S U
\end{equation}
where
$$
\overline U:=
\begin{pmatrix}
\overline U_v^1 \\
\overline U_f^1 \\
\overline U_v^2 \\
\overline U_f^2 \\
\vdots
\end{pmatrix}
$$
is the vector of own and ghost node and fracture-face unknowns on all processes. The matrix 
$S$, containing only $0$ and $1$ entries, is assembled once and for all at the beginning of the simulation. 

Finally, thanks to \eqref{eqschurcell}, the vector of own and ghost cell unknowns $\overline U_c^p$ is 
computed locally on each process $p$. 

In conclusion, the parallel implementation of the discrete hybrid dimensional Darcy flow can be summarized as:
\begin{center}
\begin{minipage}{0.82\textwidth}
\begin{algorithm}[H]
  \caption{Parallel implementation of the discrete hybrid dimensional Darcy flow}
  \label{paraimpdarcy}
  \begin{algorithmic}[1] 
    \STATE Assemble locally on each process the rectangular linear system \eqref{sysowndarcy}, 
    \STATE Compute locally on each process the Schur complement \eqref{eqschur} of \eqref{sysowndarcy}, 
    \STATE Construct the parallel linear system \eqref{parallellin} in PETSc or Trilinos,
    \STATE Solve the parallel linear system \eqref{parallellin} to obtain the solution at own nodes and fracture faces,
    \STATE Communicate the solution at ghost nodes and fracture faces from \eqref{commS},
    \STATE Compute locally on each process the solution at own and ghost cells from \eqref{eqschurcell}.
  \end{algorithmic}
\end{algorithm}
\end{minipage}
\end{center}

\subsection{Parallelization of the discrete hybrid dimensional transport model}
The parallel implementation of the transport model \eqref{modeleTracer} with an explicit upwind discretization of the fluxes consists of the following four steps. 
\begin{enumerate}
\item Compute the Darcy matrix and fracture fluxes defined by \eqref{VAGFluxMatrix} and \eqref{VAGFluxFracture}. 
\item Compute the maximum time step $\Delta t$ satisfying the CFL condition \eqref{eq_CFL} and 
set   $\Delta t^n = \Delta t$ for all 
$n=0,\cdots,N-2$, and $\Delta t^{N-1} = T-(N-1)\Delta t$ with $N= \lceil{T\over \Delta t}\rceil$. 
\item For each time step $n=0,1,...,N-2$,
\begin{enumerate}[3a.]
\item Compute $c_{\s}^{n+1}$, $c_{\sigma}^{n+1}$ and $c_{K}^{n+1}$, $\s \in  \nodes_{int} \cap \mathcal{V}^p$, $\sigma \in \mathcal{F}_\Gamma^p$, $K \in \overline{\cells}^p$ solution of the following explicit equations
\begin{equation}
\label{tracero1own}
\left\{\begin{array}{r@{\,\,}c@{\,\,}ll}
&& \dsp \phi_K {c_{K}^{n+1}- c_{K}^{n} \over \Delta t}  
+ \sum_{\s\in {\nodes}_K} H_{K,\s}(c_\D^n) + \sum_{\sigma\in {\faces}_{\Gamma,K}} H_{K,\sigma}(c_\D^n)   = 0,  
\,\, K \in \overline{\cells}^p,\\
&&  \dsp  \phi_\sigma {c_{\sigma}^{n+1}- c_{\sigma}^{n} \over \Delta t} + 
\sum_{\s \in {\nodes}_\sigma} H_{\sigma,\s}(c_\D^n) 
- \sum_{K \in \cells_\sigma} H_{K,\sigma}(c_\D^n) = 0,
\,\, \sigma \in \mathcal{F}_\Gamma^p,\\
&&\dsp \phi_\s {c_{\s}^{n+1}- c_{\s}^{n} \over \Delta t} -  \sum_{K\in \cells_\s} H_{K,\s}(c_\D^n)
- \sum_{\sigma \in \faces_{\Gamma,\s}} H_{\sigma,\s}(c_\D^n) = 0,
\,\, \s \in \nodes_{int} \cap \mathcal{V}^p,\\
&& c_\s = \bar c_\s, \,\, \s \in  \nodes_{ext} \cap \mathcal{V}^p.
\end{array}\right.
\end{equation}
\item Get the node and fracture-face ghost unknowns $c_{\s}^{n+1}$, $c_{\sigma}^{n+1}$, $\s \in  \nodes_{int} \cap (\overline{\mathcal{V}}^p \backslash \mathcal{V}^p)$, $\sigma \in  \overline{\mathcal{F}}_\Gamma^p \backslash \mathcal{F}_\Gamma^p$ using 
the PETSc or Trilinos matrix vector product with the matrix $S$ defined in \eqref{commS}, as was done for the Darcy flow solution $U$: 
$$
\begin{pmatrix}
\overline C_v^1 \\
\overline C_f^1 \\
\overline C_v^2 \\
\overline C_f^2 \\
\vdots
\end{pmatrix}
= S
\begin{pmatrix}
C_v^1 \\
C_f^1 \\
C_v^2 \\
C_f^2 \\
\vdots
\end{pmatrix}
$$
\end{enumerate}
\end{enumerate}
Thanks to our mesh decomposition, step 1 and step 3a are performed locally on each process. 
For step 2, the maximum time step $\Delta t^p$ is computed locally on each process $p$, then the time step 
$\Delta t$ is obtained using the MPI reduce operation.

\section{Numerical experiments}
\label{sec_tests}

All the numerical tests have been implemented on the Cicada cluster of the University Nice Sophia Antipolis 
consisting of 72 nodes (16 cores/node, Intel Sandy Bridge E5-2670, 64GB/node). 
We always fix 1 core per process and 16 processes per node. 
The communications are handled by OpenMPI 1.8.2 (GCC 4.9) and PETSc 3.5.3. 

The first two test cases are designed in order to validate the Darcy fluxes and the convergence 
of the transport model discretization on two analytical solutions including one fracture for the first test case and four intersecting fractures for the second test case. 
In the remaining test cases, the parallel scalability of our Darcy flow and transport solvers   
is assessed with different types of fracture networks and meshes and different matrix fracture permeability ratios. 
In particular, the last test case applies our algorithm to a complex fracture network with hundreds of fractures.

\subsection{Numerical convergence for an analytical solution with one fracture}
\begin{figure}[!htbp]
\centering
\includegraphics[width=0.35\textwidth]{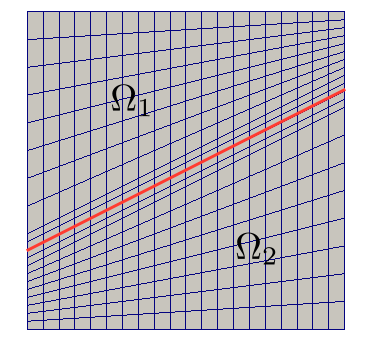}
\hspace{0.05\textwidth}
\includegraphics[width=0.4\textwidth]{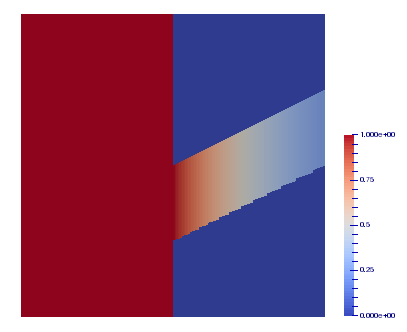}
\caption{Left: example of mesh with $n_x=20$ where the red line is the fracture. 
Right: analytical solution of \eqref{test_traceur} at time $t_f = 0.5$.}
\label{cpgmeshsol}
\end{figure}

\begin{figure}[!htbp]
\centering
\includegraphics[width=0.45\textwidth]{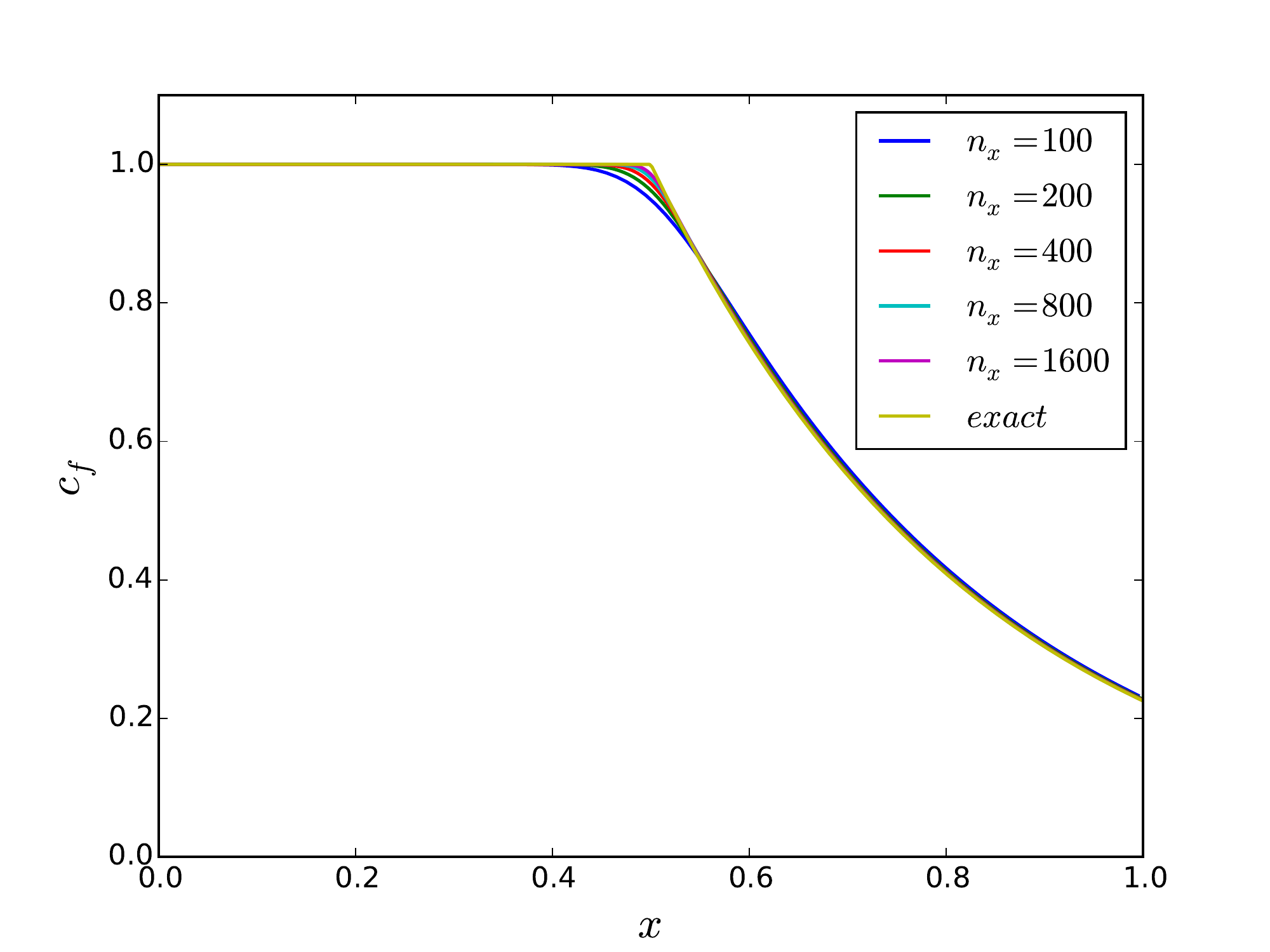}
\hspace{0.05\textwidth}
\includegraphics[width=0.45\textwidth]{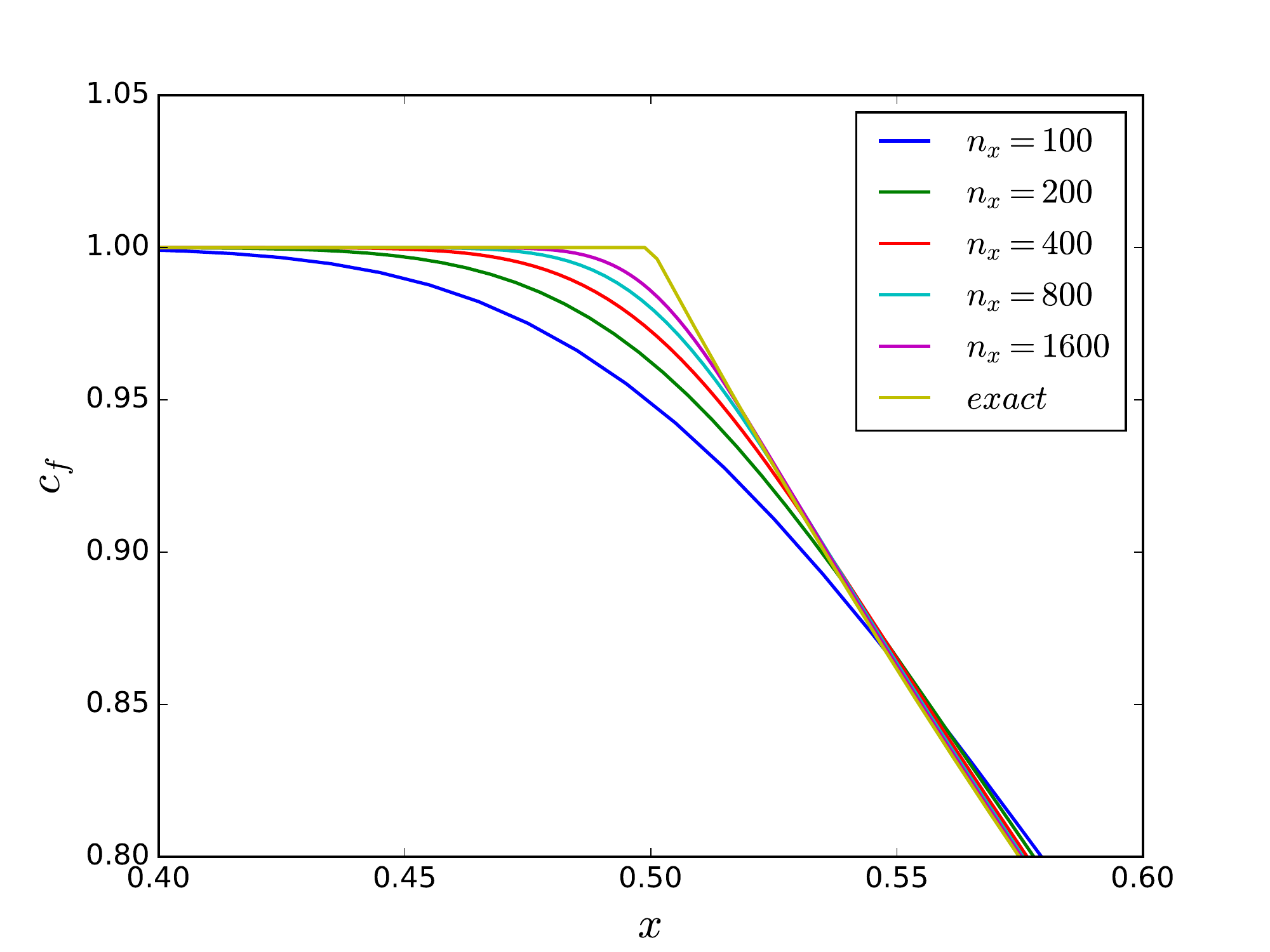}
\caption{Left: analytical solution and numerical solutions along the fracture at time $t_f$ with $n_x=100,200,400,800,1600$. Right: zoom view of left figure around $x=0.5$.}
\label{solonfrac}
\end{figure}

\begin{figure}[!htbp]
\centering
\includegraphics[width=0.5\textwidth]{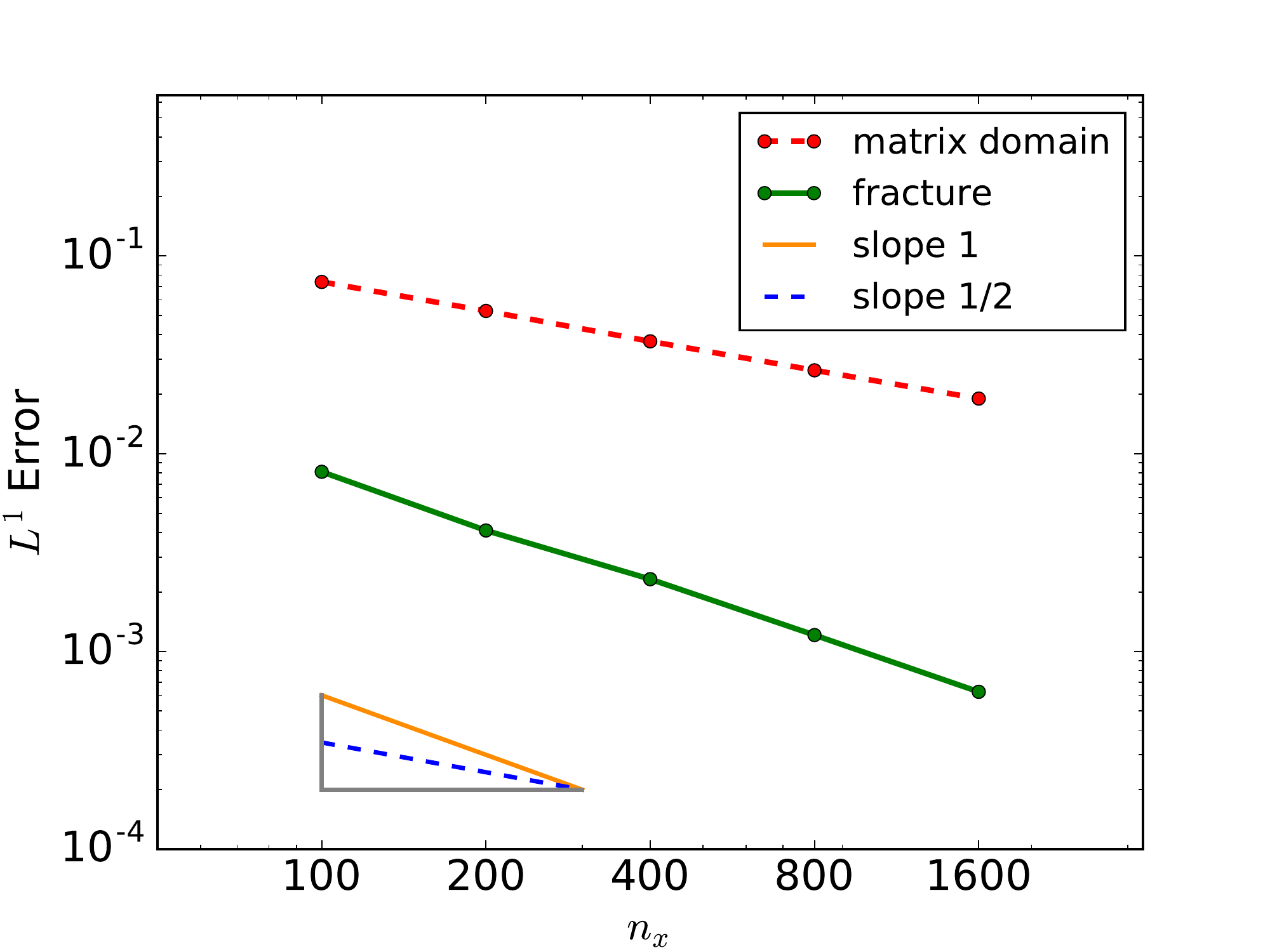}
\caption{Relative $L^1$ errors and references in the matrix domain and in the fracture at time $t_f$ 
between the analytical solution and the numerical solutions as a function of the grid size 
 $n_x=100,200,400,800,1600$.}
\label{erro1o2frac}
\end{figure}

Let us set $\Omega = (0,1)^2$, and denote by $(x,y)$ the Cartesian coordinates of $\x$.
We then define $\x_1 = (0,{1\over 4})$, $\theta\in (0,\mbox{arctan}({3\over 4}))$, 
$\x_2 = (1,{1\over 4}+\tan(\theta))$. Let $\Omega_1 = \{(x,y)\in\Omega \,|\, y > {1\over 4} + x\tan(\theta)\}$, 
and  $\Omega_2 = \Omega\setminus \overline\Omega_1$. 
We consider a single fracture defined by 
$\Gamma = (\x_1,\x_2)=\partial \Omega_1\cap \partial\Omega_2$ with tangential permeability 
$\Lambda_f > 0$, and width $d_f>0$. The matrix permeability is isotropic and set to $\Lambda_m = 1$, the  matrix and 
fracture porosities are set to $\phi_m=\phi_f = 1$, and  the fluid viscosity is set to $\mu=1$. 
The pressure solution is fixed to $u(x,y) = 1-x$. 
In this case, the transport model \eqref{modeleTracer} reduces to the following system of equations which specifies our choice 
of the boundary and initial conditions: 
\begin{eqnarray}
\label{test_traceur}
\left\{\begin{array}{r@{\,\,}c@{\,\,}ll}
\partial_t c_{m,\alpha}(x,y,t) + \partial_x c_{m,\alpha}(x,y,t) &=& 0 &\mbox{ on } \Omega_\alpha\times (0,T), \,\alpha=1,2,\\
c_{m,\alpha}(x,y,0) &=& 0 &\mbox{ on } \Omega_\alpha, \,\alpha=1,2,\\  
c_{m,1}(0,y,t) &=& 1 &\mbox{ on } ({1\over 4},1)\times (0,T),\\
c_{m,2}(0,y,t) &=& 1 &\mbox{ on } (0,{1\over 4})\times (0,T),\\
c_{m,2}(x,{1\over 4} + x\tan(\theta),t) &=& c_f(x,t) &\mbox{ on } (0,1)\times (0,T),\\
{\cal L} c_f(x,t) &=& 
\beta c_{m,1}(x,{1\over 4} + x\tan(\theta),t)
&\mbox{ on } (0,1)\times (0,T), \\
c_f(0,t) &=& 1 &\mbox{ on } (0,T),\\
c_f(x,0) &=& 0 &\mbox{ on } (0,1),
\end{array}\right.
\end{eqnarray}
where  ${\cal L} = \partial_t + k \partial_x + \beta$ with  $\beta = {\sin(\theta)\over d_f}$ and 
$k=\Lambda_f \cos^2(\theta)$. It is assumed that $k>1$.  
This system can be integrated along the characteristics of the matrix and fracture velocity fields leading 
to the following analytical solution: 
$$
c_{m,1}(x,y,t) = 
\left\{\begin{array}{r@{\,\,}c@{\,\,}ll}
&0 &\mbox{if}& t < x,\\
&1 &\mbox{if}& t > x,
\end{array}\right.
$$
$$
c_f(x,t) = 
\left\{\begin{array}{r@{\,\,}c@{\,\,}ll}
&0 &\mbox{if} & t < {x\over k},\\
&e^{-{\beta\over k-1}(x-t)} &\mbox{if}& {x\over k} < t < x,\\
&1 &\mbox{if} & t > x,
\end{array}\right.
$$
$$
c_{m,2}(x,y,t) = 
\left\{\begin{array}{r@{\,\,}c@{\,\,}l}
&&\mbox{ if } y\in(0,{1\over 4})
\left\{\begin{array}{r@{\,\,}c@{\,\,}ll}
&0 &\mbox{if}& t < x,\\
&1 &\mbox{if}& t > x,
\end{array}\right.\\
&&\mbox{ if } y\in({1\over 4}, {1\over 4} + \tan(\theta))
\left\{\begin{array}{r@{\,\,}c@{\,\,}ll}
&0 &\mbox{if}& t < x-{4y-1\over 4\tan(\theta)},\\
&c_f({4y-1\over 4\tan(\theta)},t+{4y-1\over 4\tan(\theta)}-x) &\mbox{if}& 
t > x-{4y-1\over 4\tan(\theta)}. 
\end{array}\right. 
\end{array}\right. 
$$
In the following numerical experiments the parameters are set to $\tan(\theta) = {1\over 2}$, $\Lambda_f = 20$ and $d_f = 0.01$.  
The mesh is a topologically Cartesian $n_x\times n_x$ grid. Figure \ref{cpgmeshsol} shows an example of the mesh with $n_x=20$ as well 
as the analytical solution in the matrix 
obtained at time $t_f=0.5$ chosen as the final time of the simulation. 
The time step is defined by the maximum time step allowed by the CFL condition \eqref{eq_CFL}. 
Figure \ref{erro1o2frac} exhibits the convergence of the relative $L^1$ errors between the analytical solution and 
the numerical solution at time $t_f$ both 
in the matrix domain and in the fracture as a function of the grid size 
$n_x=100,200,400,800,1600$. Figure \ref{solonfrac} shows the analytical solution and the numerical solutions obtained 
at time $t_f$ along the fracture. In both cases, we can observe the expected convergence of the numerical solution 
to the analytical solution with a higher order of convergence in the fracture due to the fact that at time $t_f$ the 
analytical solution in the fracture is continuous as exhibited in Figure \ref{solonfrac}. 

\subsection{Numerical convergence for an analytical solution with four intersecting fractures}

Let $\Omega = (0,1)^2$, $\x_1 = (0,{1\over 4})$, $\theta_1\in (0,\mbox{arctan}({3\over 4}))$, 
$\x_2 = (1,{1\over 4}+\tan(\theta_1))$, 
$\x_3 = ({3\over 4},0)$, $\x_4 = ({3\over 4}-\tan(\theta_2),1)$, and the 
intersection of $\x_1\x_2$ and $\x_3\x_4$ equal to 
$$
\x_0 = (x_0,y_0) = {1\over 4(1+\tan(\theta_1)\tan(\theta_2))}(3-\tan(\theta_2), 1+3\tan(\theta_1)).
$$ 
It is assumed that $\theta_1,\theta_2\in (0,\mbox{arctan}({3\over 4}))$. 

We consider the four fractures
$\Gamma_1 = (\x_1,\x_0)$, $\Gamma_2 = (\x_0,\x_2)$, $\Gamma_3 = (\x_3,\x_0)$, $\Gamma_4 = (\x_0,\x_4)$, with tangential 
permeabilities $\Lambda_{f,1}=\Lambda_{f,2} > 0$, and $\Lambda_{f,3}=\Lambda_{f,4} > 0$, and with widths $d_{f,1} = d_{f,2} >0$, and $d_{f,3} = d_{f,4} >0$.  It is assumed that $\Lambda_m=1$. 

The fractures partition the domain $\Omega$ in the following four subdomains 
$$
\Omega_1 = \{\x=(x,y)\in \Omega \,|\, y > {1\over 4} + x\tan(\theta_1), x <  {3\over 4} - y\tan(\theta_2)\},
$$
$$
\Omega_2 = \{\x=(x,y)\in \Omega \,|\, y > {1\over 4} + x\tan(\theta_1), x >  {3\over 4} - y\tan(\theta_2)\},
$$
$$
\Omega_3 = \{\x=(x,y)\in \Omega \,|\, y < {1\over 4} + x\tan(\theta_1), x <  {3\over 4} - y\tan(\theta_2)\},
$$
$$
\Omega_4 = \{\x=(x,y)\in \Omega \,|\, y < {1\over 4} + x\tan(\theta_1), x >  {3\over 4} - y\tan(\theta_2)\}. 
$$
Let us set $\beta_1 = {\sin(\theta_1)\over d_{f,1}}$, 
$k_1=\Lambda_{f,1} \cos^2(\theta_1)$, $\beta_2 = {\cos(\theta_2)\over d_{f,3}}$, 
$k_2=\Lambda_{f,3} \cos(\theta_2)\sin(\theta_2)$, 
$r = {\Lambda_{f,3}d_{f,3} \sin(\theta_2) \over \Lambda_{f,1}d_{f,1} \cos(\theta_1)}$. 
It is assumed that $k_1>1$ and $k_2 \tan(\theta_2) > 1$.  The  matrix and 
fracture porosities are set to $\phi_m=\phi_f = 1$ and  the fluid viscosity is set to $\mu=1$. \\

The pressure solution is set to $u(x,y) = 1-x$. In that case, the transport model \eqref{modeleTracer} reduces to the 
following system of equations which specifies our choice of the boundary and initial conditions: 
find $c_{m,\alpha}(x,y,t)$, $\alpha=1,\cdots,4$, $c_{f,1}(x,t)$,  $c_{f,2}(x,t)$,  $c_{f,3}(y,t)$,  $c_{f,4}(y,t)$, 
and $c_0(t)$ such that   
\begin{eqnarray}
\label{test_traceur4frac}
\left\{\begin{array}{r@{\,\,}c@{\,\,}ll}
\partial_t c_{m,\alpha}(x,y,t) + \partial_x c_{m,\alpha}(x,y,t) &=& 0 &\mbox{ on } \Omega_\alpha\times (0,T), 
\,\alpha=1,\cdots,4,\\
c_{m,\alpha}(x,y,0) &=& 0 &\mbox{ on } \Omega_\alpha, \,\alpha=1,\cdots,4,\\  
c_{m,1}(0,y,t) &=& 0 &\mbox{ on } ({1\over 4},1)\times (0,T),\\
c_{m,3}(0,y,t) &=& 0 &\mbox{ on } (0,{1\over 4})\times (0,T),\\
c_{m,2}({3\over 4} - y\tan(\theta_2),y,t) &=& c_{f,4}(y,t) &\mbox{ on } (y_0,1)\times (0,T),\\
c_{m,4}({3\over 4} - y\tan(\theta_2),y,t) &=& c_{f,3}(y,t) &\mbox{ on } (0,y_0)\times (0,T),\\
c_{m,3}(x,{1\over 4} + x\tan(\theta_1),t) &=& c_{f,1}(x,t) &\mbox{ on } (0,x_0)\times (0,T),\\
c_{m,4}(x,{1\over 4} + x\tan(\theta_1),t) &=& c_{f,2}(x,t) &\mbox{ on } (x_0,1)\times (0,T),\\
{\cal L}_1 c_{f,1}(x,t)  -  \beta_1 c_{m,1}(x,{1\over 4} + x\tan(\theta_1),t)  &=& 0
&\mbox{ on } (0,x_0)\times (0,T), \\
{\cal L}_1 c_{f,2}(x,t)  - \beta_1 c_{m,2}(x,{1\over 4} + x\tan(\theta_1),t) &=& 0
&\mbox{ on } (x_0,1)\times (0,T), \\
{\cal L}_2 c_{f,3}(y,t) - \beta_2 c_{m,3}({3\over 4} - y\tan(\theta_2),y,t) &=& 0
&\mbox{ on } (0,y_0)\times (0,T), \\
{\cal L}_2 c_{f,4}(y,t)  - \beta_2 c_{m,1}({3\over 4} - y\tan(\theta_2),y,t) &=& 0
&\mbox{ on } (y_0,1)\times (0,T), \\
c_{f,2}(x_0,t) = c_{f,3}(y_0,t) &=& c_0(t) &\mbox{ on } (0,T),\\
(r+1) c_0(t)   - c_{f,1}(x_0,t) - r  c_{f,4}(y_0,t) &=& 0 &\mbox{ on } (0,T),\\
c_{f,1}(0,t)  = c_{f,4}(1,t)  &=& 1 &\mbox{ on } (0,T),\\
c_{f,1}(x,0) &=& 0 &\mbox{ on } (0,x_0),\\
c_{f,2}(x,0) &=& 0 &\mbox{ on } (x_0,1),\\
c_{f,3}(y,0) &=& 0 &\mbox{ on } (0,y_0),\\
c_{f,4}(y,0) &=& 0 &\mbox{ on } (y_0,1),  
\end{array}\right.
\end{eqnarray}
 where ${\cal L}_1 = \partial_t + k_1 \partial_x + \beta_1$ and ${\cal L}_2 = \partial_t - k_2 \partial_y + \beta_2$.

This system can also be integrated analytically along the characteristics of the matrix and fracture velocity fields, 
but it leads to complex computations. It is much easier to obtain the stationary solution of this system which 
is defined in the fractures by 
$$
c_{f,1}(x) = e^{-{\beta_1\over k_1} x},
$$
$$
c_{f,4}(y) = e^{-{\beta_2\over k_2 }(1-y) },
$$
$$
c_0 = { e^{-{\beta_1\over k_1} x_0} + r e^{-{\beta_2\over k_2 }(1-y_0)} \over r+1},
$$
$$
c_{f,2}(x) = e^{-{\beta_1\over k_1} x} \( c_0  e^{{\beta_1\over k_1} x_0} 
+ {\beta_1\over k_1 r_1} 
\( e^{ (r_1 x-{3\beta_2\over 4 k_2}) } 
-  e^{ (r_1 x_0-{3\beta_2\over 4 k_2}) } 
\) \),
$$
$$
c_{f,3}(y) = 
\left\{\begin{array}{r@{\,\,}c@{\,\,}ll}
& e^{{\beta_2\over k_2} y} \( c_0  e^{-{\beta_2\over k_2} y_0} 
+ {\beta_2\over k_2 r_2} 
\( 
e^{ ( - {r_2\over 4}-{\beta_1\over 4 k_1\tan(\theta_1)}) } 
-e^{ ( -r_2 y_0-{\beta_1\over 4 k_1\tan(\theta_1)}) } 
\) \),&\mbox{if}& y < {1\over 4},\\
& e^{{\beta_2\over k_2} y} \( c_0  e^{-{\beta_2\over k_2} y_0} 
+ {\beta_2\over k_2 r_2} 
\( 
e^{ ( -r_2 y-{\beta_1\over 4 k_1\tan(\theta_1)}) } 
-e^{ ( -r_2 y_0-{\beta_1\over 4 k_1\tan(\theta_1)}) } 
\) \)&\mbox{if}& y > {1\over 4},
\end{array}\right.
$$
with $r_1 = {\beta_1\over k_1}+{\beta_2\over k_2}\tan(\theta_1)$ and $r_2 = {\beta_2\over k_2}+{\beta_1\over k_1\tan(\theta_1)}$, 
and in the matrix by 
$$
c_{m,1}(x,y) = 0, 
$$
$$
c_{m,2}(x,y) = c_{f,4}(y),
$$
$$
c_{m,3}(x,y) = 
\left\{\begin{array}{r@{\,\,}c@{\,\,}ll}
&0 &\mbox{if}& y < {1\over 4},\\
& c_{f,1}\({y-{1\over 4}\over \tan(\theta_1)}\) &\mbox{if}& y > {1\over 4},
\end{array}\right.
$$

$$
c_{m,4}(x,y) = 
\left\{\begin{array}{r@{\,\,}c@{\,\,}ll}
&c_{f,3}(y) &\mbox{if}& y < y_0,\\
& c_{f,2}\({y-{1\over 4}\over \tan(\theta_1)}\) &\mbox{if}& y > y_0. 
\end{array}\right.
$$

\begin{figure}[!htbp]
\centering
\includegraphics[height=5.5cm]{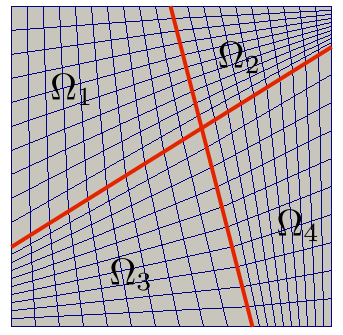}
\hspace{0.05\textwidth}
\includegraphics[height=5.5cm]{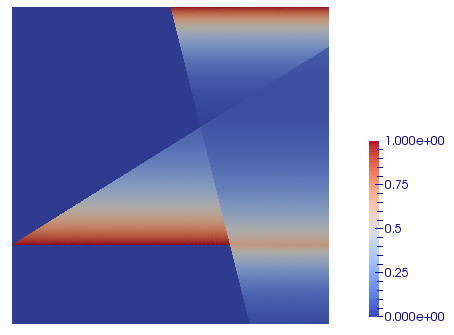}
\caption{Left: example of mesh with $n_x=20$ where the red lines account for the four fractures. 
Right: stationary analytical solution of \eqref{test_traceur4frac}.}
\label{cpg4fracmeshsol}
\end{figure}

In the following numerical experiments the parameters 
are set to $\tan(\theta_1) = {5\over 8}$, $\tan(\theta_2) = {1\over 4}$, 
$\Lambda_{f,1} = 200$, $\Lambda_{f,3} = 400$, and $d_{f,1} = d_{f,3}  = 0.01$.  
The mesh is, as for the previous test case, a topologically Cartesian $n_x\times n_x$ grid exhibited in
Figure \ref{cpg4fracmeshsol} for $n_x = 20$. Figure \ref{cpg4fracmeshsol} also shows  
the stationary analytical solution in the matrix. 
The time step is again defined by the maximum time step allowed by the CFL condition \eqref{eq_CFL} 
and the simulation time is chosen large enough to obtain the numerical stationary solution.  

\begin{figure}[!htbp]
\centering
\includegraphics[width=0.45\textwidth]{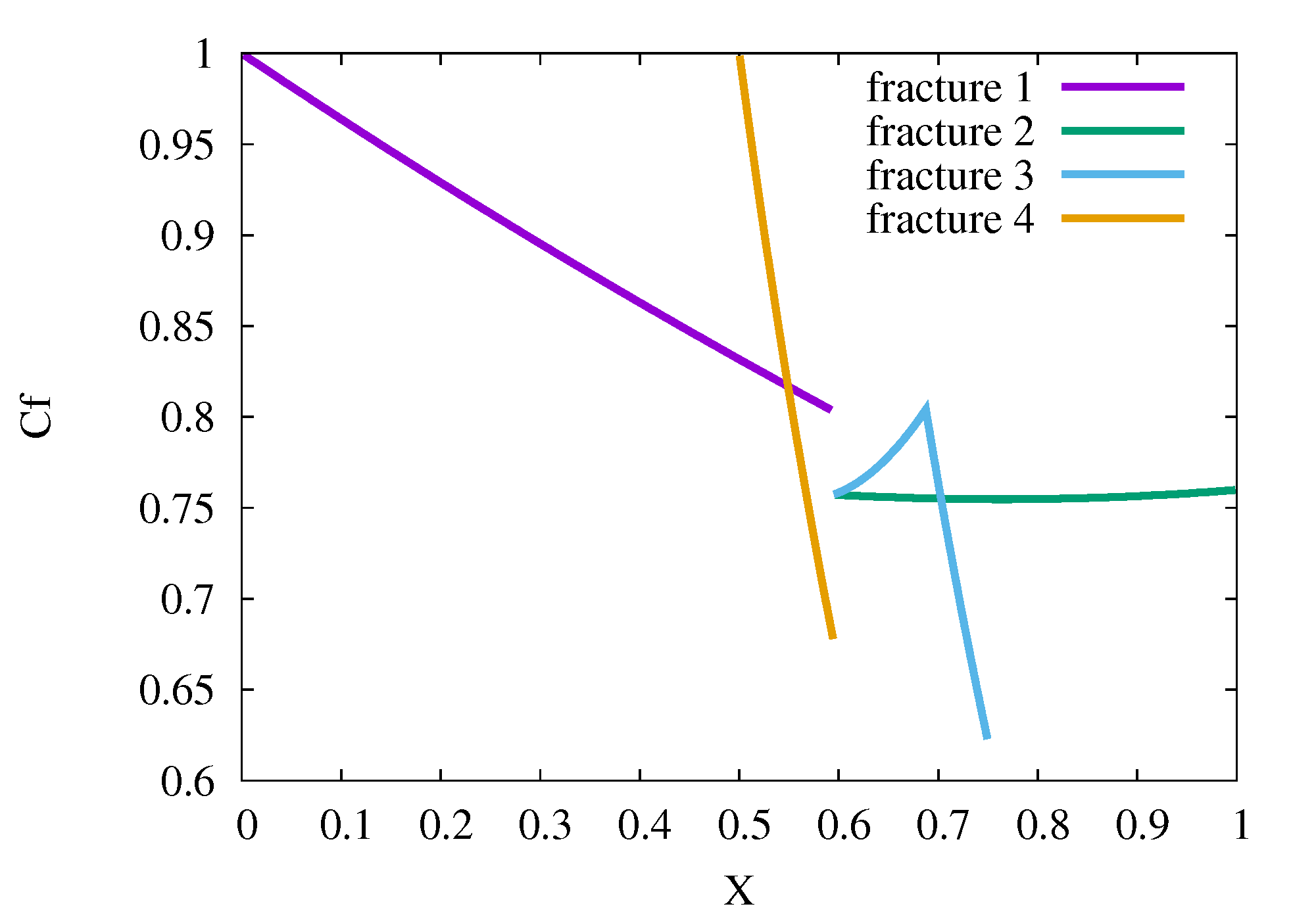}
\includegraphics[width=0.45\textwidth]{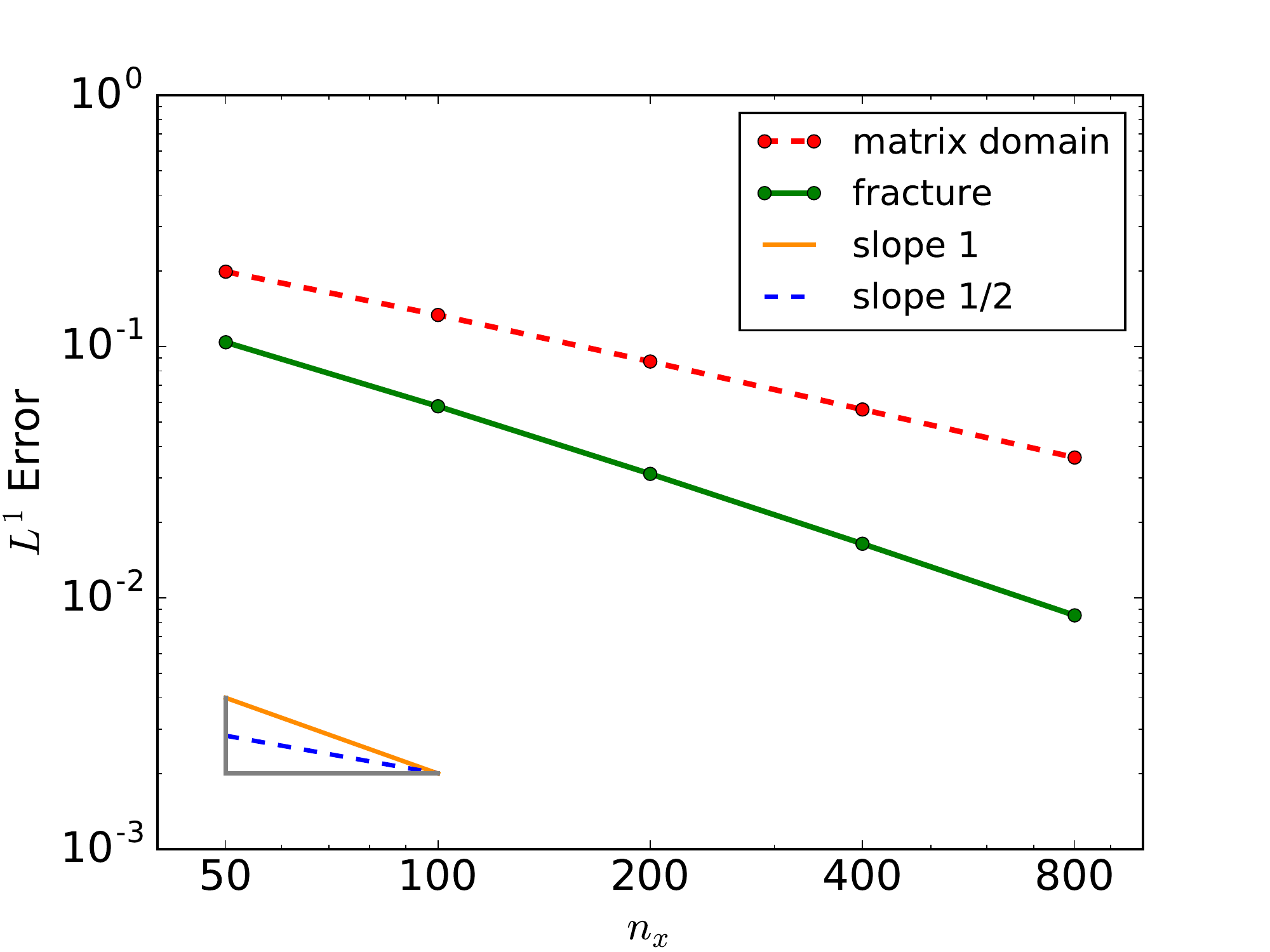}
\caption{Left: stationary analytical solution in the four fractures as a function of the $x$ coordinate. 
Right: relative $L^1$ errors and references in the matrix domain and in the fracture network 
between the stationary analytical and numerical solutions as a function of the grid size 
 $n_x=100,200,400,800$.}
\label{erro1o2frac4}
\end{figure}

Figure \ref{erro1o2frac4} exhibits the convergence of the relative $L^1$ errors between the stationary analytical and 
the numerical solutions  both 
in the matrix domain and in the fracture as a function of the grid size 
$n_x=100,200,400,800$. We can again observe the expected convergence of the numerical solution 
to the analytical solution with a higher order of convergence in the fracture network due to the fact that 
the solution is continuous on each individual fracture as exhibited in Figure \ref{erro1o2frac4}. This property 
is always true when looking at the solution at the matrix time scale and could be exploited in a future work 
by using an implicit time integration in the fracture coupled to an explicit time integration in the matrix 
domain with a higher order discretization in space in the spirit of what has been done in \cite{HF08}.

\subsection{Fracture network with hexahedral meshes}

\begin{figure}[!htbp]
\centering
\includegraphics[width=0.38\textwidth]{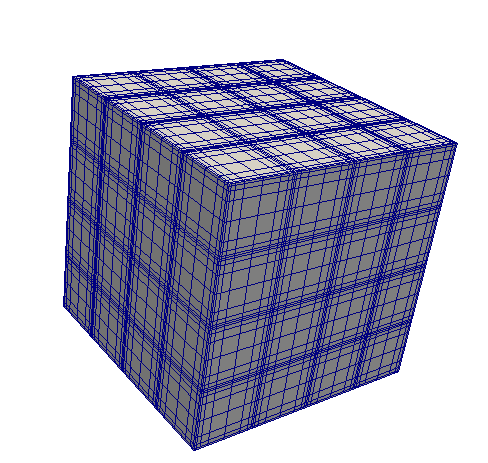}
\hspace{0.05\textwidth}
\includegraphics[width=0.38\textwidth]{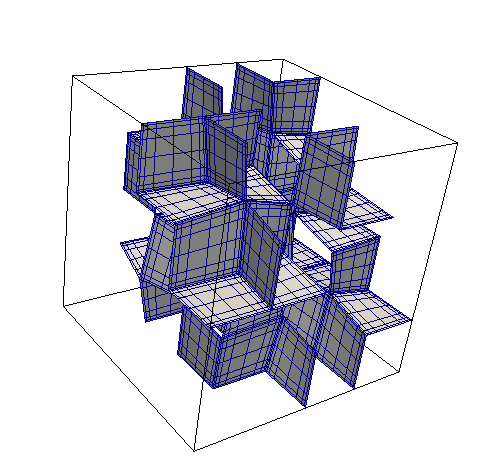}
\caption{Hexahedral mesh in the matrix domain (left) and in the fracture network (right) with $n_x=32$.}
\label{meshcpgfrac}
\end{figure}

\begin{figure}[!htbp]
          \centering
          \addtocounter{subfigure}{2} 
          \begin{subfigure}{0.48\textwidth}
           \centering
            \includegraphics[width=\textwidth]{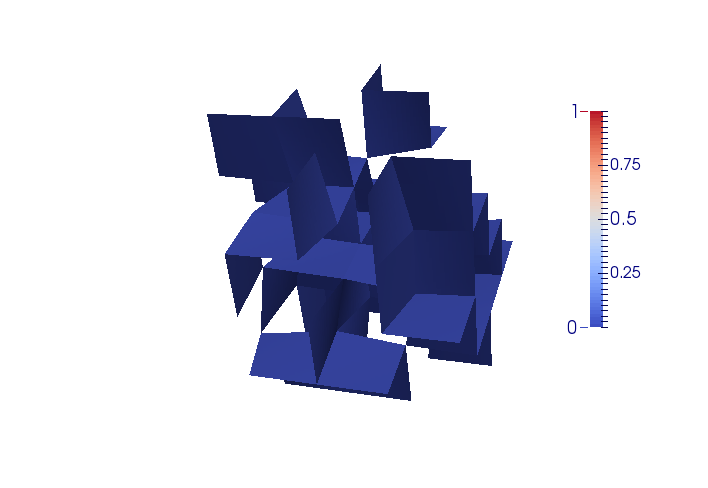}
            \caption{$t=0.0$}
          \end{subfigure}
          \begin{subfigure}{0.48\textwidth}
           \centering
            \includegraphics[width=\textwidth]{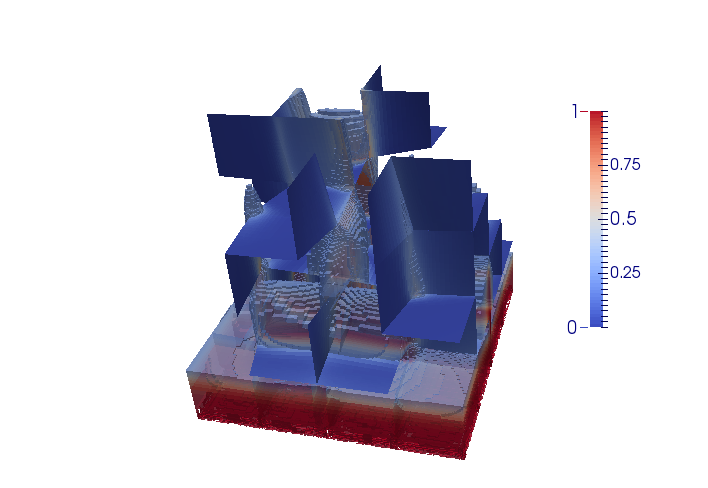}
            \caption{$t=0.2$}
          \end{subfigure}
          \begin{subfigure}{0.48\textwidth}
           \centering
            \includegraphics[width=\textwidth]{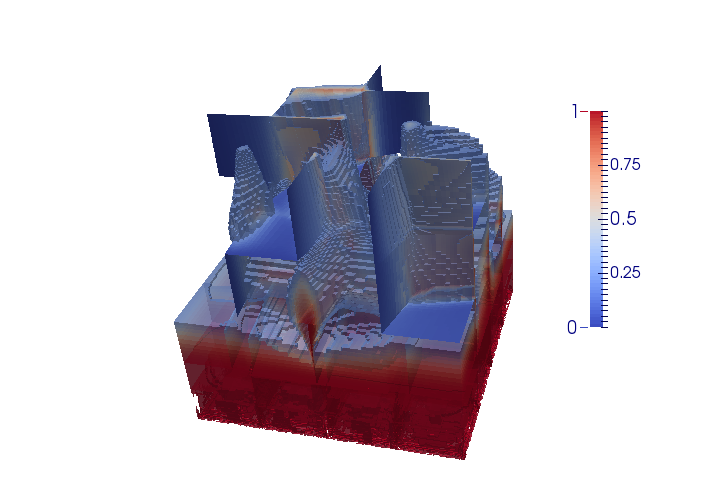}
            \caption{$t=0.4$}
          \end{subfigure}
          \begin{subfigure}{0.48\textwidth}
           \centering
            \includegraphics[width=\textwidth]{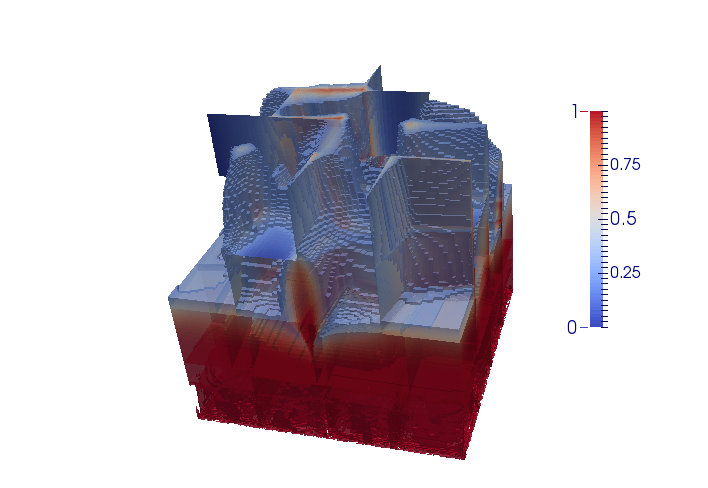}
            \caption{$t=0.5$}
          \end{subfigure}
          \caption{Concentration in the matrix domain and in the fracture network obtained 
                  at different times with the mesh $n_x = 128$. 
                   A threshold concentration of 0.2 is used in the matrix domain.}
          \label{simucpgfrac}
        \end{figure}

The objective of this subsection and of the next subsection is 
to investigate the parallel scalability of the algorithms described in Section \ref{sec_parallel}. 
In this subsection we consider a topologically Cartesian mesh of size $n_x\times n_x\times n_x$ of the cubic domain $(0,1)^3$ 
exhibited in Figure \ref{meshcpgfrac} for $n_x = 32$. 
The mesh is exponentially refined at the interface between the 
matrix and the fracture network exhibited in Figure \ref{simucpgfrac}. 
The permeabilities are isotropic and set to $\Lambda_f=20$ in the fracture network and to 
$\Lambda_m=1$ in the matrix. The porosities are set to $\phi_m = \phi_f = 1$ and the fluid viscosity 
is set to $\mu = 1$. The initial concentration is set to $0$ both in the matrix domain 
and in the fracture network and a value of $1$ of the concentration is injected on the bottom boundaries of the matrix 
and of the fracture network. The pressure is fixed to $u=1$ and $\gamma u=1$ on the bottom boundary and 
to $u=0$ and $\gamma u = 0$ on the top boundary. The remaining lateral boundaries are considered impervious. 
Figure \ref{simucpgfrac} exhibits the tracer concentrations obtained with the mesh $n_x=128$ 
at times $t=0$, $t=0.2$, $t=0.4$ and at the final simulation time $t_f = 0.5$.  

Table \ref{solver_cpg_nb} presents the numbers of linear solver iterations for the stationary pressure solution for a number of MPI processes 
ranging from $N_p = 2$ to $N_p = 512$ and with the mesh size $n_x = 128$ 
corresponding to roughly $2.1\times 10^6$ cells, $2.1\times 10^6$ nodes and $5.2\times 10^4$ fractures faces.
Both the GMRES and BiCGStab linear solvers from the PETSc library 
are tested combined with either 
the Boomer AMG preconditioner from the Hypre library \cite{hyprebib}, 
the Aggregation AMG preconditioner 
from the Trilinos library \cite{citetrilinos} 
or the block Jacobi ILU(0) preconditioner from the Euclid library. No restart is used 
for the GMRES linear solver. 
Table \ref{solver_cpg_time} shows the corresponding computation times both for the 
setup phase of the preconditioner and for the solve phase of the linear solver. 

 According to these tables, the GMRES and the BiCGStab linear solvers combined with 
the Boomer AMG preconditioner are good choices for a number of processes 
$N_p \leqslant 128$, while the BiCGStab linear solver combined with the block Jacobi ILU(0) 
preconditioner is more efficient for this test case for $N_p=256$ and $N_p=512$. 
This was expected since the Boomer AMG preconditioner requires a sufficiently large number 
of unknowns per core to maintain a good parallel scalability due to the high level 
of communications in particular in the setup phase of the preconditioner.  
For this linear system, the number of unknowns per MPI process is roughly 
$4100$ for $N_p = 512$  which is too small for this type of preconditioner while 
the block Jacobi preconditioner still maintains a good parallel scalability for such a number  
of unknowns per MPI process. On the other hand, as shown in Table \ref{solver_cpg_perm}, 
Boomer AMG exhibits an optimal scalability 
while ILU(0) is not scalable 
in terms of iteration count with respect to the ratio 
${\Lambda_f \over \Lambda_m}$ between the fracture and matrix permeabilities. 
The same remark also holds in terms of scalability with respect to the mesh size which 
means that the ILU(0) preconditioner is only advantageous for small size and mildly heterogeneous 
problems.  

Tables \ref{solver_cpg_nb} and \ref{solver_cpg_time} also clearly show  
that the BiCGStab linear solver  
outperforms the GMRES linear solver for 
cases requiring a large number of iterations due to the fact that the cost 
of the orthogonalization procedure increases with the dimension of the Krylov subspace. 
The Aggregation AMG preconditioner yields a larger number of iterations 
compared with the Boomer AMG preconditioner but has a much lower setup time resulting for 
this test case in a total lower CPU time.  However, this implementation of the 
Aggregation AMG preconditioner seems to lack robustness with respect to the matrix fracture 
permeability ratio as exhibited in Table \ref{solver_cpg_perm}.

\begin{table}[!htbp]
\caption{Number of linear solver iterations vs. the number of MPI processes obtained with  different linear solvers and 
preconditioners for the mesh size $n_x=128$.}
\label{solver_cpg_nb}
\centering
\small
\begin{tabular}{c|c|c|c|c|c|c|c|c|c}
\hline
$N_p$ & 2 & 4 & 8 & 16 & 32 & 64 & 128 & 256 & 512 \\
\hline
GMRES + Boomer AMG & 15 & 15  & 15 & 15   & 15   & 16   & 15   & 15   & 15  \\
\hline
GMRES + Aggregation AMG & 59 & 78 & 65 & 39 & 65 & 54 & 73 & 53 & 62 \\
\hline
GMRES + ILU(0) & 751 & 707 & 655 & 644 & 648 & 634 & 633 & 624 & 613\\ 
\hline
BiCGStab + Boomer AMG & 9 & 9 & 9 & 9 & 9 & 10 & 9 & 9 & 10 \\
\hline
BiCGStab + ILU(0) & 508  & 476 & 484 & 503 & 473 & 513 & 491 & 487 & 484 \\
\hline
\end{tabular}
\end{table}

\begin{table}[!htbp]
\caption{Linear solver setup phase and solve phase computation times vs. the number of MPI processes 
obtained with  different linear solvers and 
preconditioners for the mesh size $n_x=128$.}
\label{solver_cpg_time}
\centering
\small
\begin{tabular}{c|c|c|c|c|c|c|c|c|c|c}
\hline
$N_p$ & & 2 & 4 & 8 & 16 & 32 & 64 & 128 & 256 & 512 \\
\hline
\multirow{2}{3.1cm}{GMRES\\ Boomer AMG} 
& Setup & 34.1 & 20.1 & 16.3 & 11.7 & 11.3 & 11.9 & 12.1 & 19.2 & 29.3 \\
& Solve & 26.3 & 15.6 & 14.8 & 7.2  & 5.2  & 3.8  & 2.5  & 5.2  & 9.6 \\
\hline
\multirow{2}{3.1cm}{GMRES \\ Aggregation AMG} 
& Setup & 4.7 & 1.9 & 1.6 & 1.5 & 2.3 & 2.9 & 4.4 & 6.6 & 11.3 \\
& Solve & 45.1 & 20.9 & 17.0 & 9.7 & 5.2 & 2.5 & 2.3 & 1.5 & 3.2 \\
\hline
\multirow{2}{3.1cm}{GMRES \\ ILU(0)} 
& Setup & 16.9 & 21.3 & 16.3 & 23.2 & 14.6 & 11.0 & 9.7 & 6.0 & 4.8 \\
& Solve & 672.3 & 590.9 & 281.6 & 163.9 & 71.4 & 30.7 & 16.7 & 8.3 & 4.0 \\
\hline
\multirow{2}{3.1cm}{BiCGStab \\ Boomer AMG} 
& Setup & 38.0 & 23.3 & 15.0 & 10.3 & 9.1 & 9.4 & 12.8 & 14.8 & 23.8 \\
& Solve & 37.1 & 21.3 & 11.5 & 7.4 & 4.1 & 2.9 & 2.5 & 4.4 & 10.0 \\
\hline
\multirow{2}{3.1cm}{BiCGStab \\ ILU(0)} 
& Setup & 18.9 & 19.9 & 16.5 & 22.1 & 14.3 & 12.4 & 9.4 & 5.8 & 3.9 \\
& Solve & 179.4 & 111.7 & 86.0 & 59.9 & 27.9 & 15.4  & 8.0 & 4.2 & 2.2 \\
\hline
\end{tabular}
\end{table}

\begin{figure}[!htbp]
\centering
\includegraphics[width=0.45\textwidth]{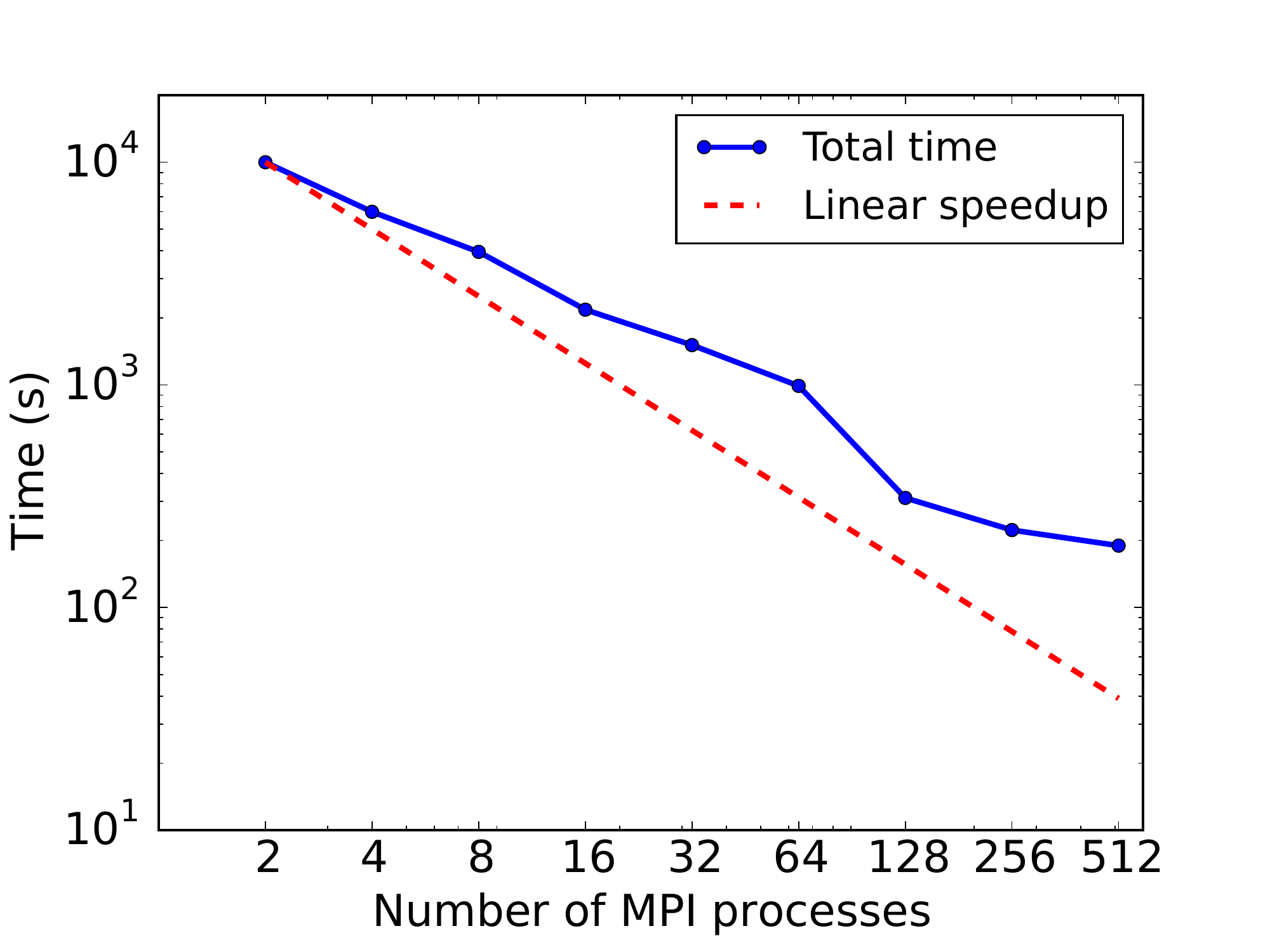}
\hspace{0.05\textwidth}
\includegraphics[width=0.45\textwidth]{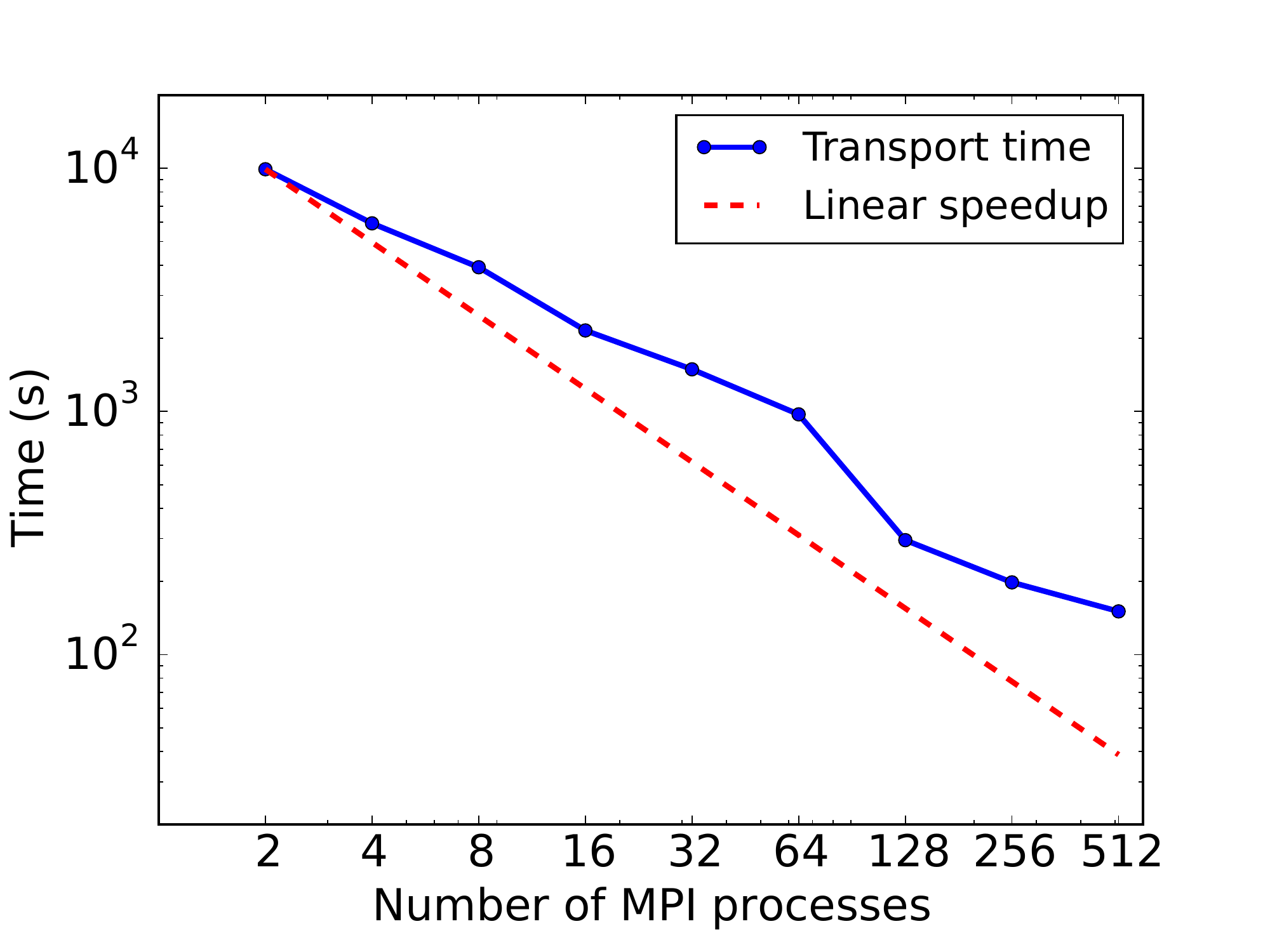}
\caption{Total computation time (left) and computation time for the transport model (right) 
vs. the number of MPI processes for the mesh size $n_x=128$.}
\label{timecpgfrac}
\end{figure}

Next, Figure \ref{timecpgfrac} plots the total (Darcy flow and transport models) computation time 
and the computation time for the transport model only as a function of the number of processes. 
In these runs the GMRES linear solver is used combined with the 
Boomer AMG preconditioner for $N_p \leqslant 128$ and with the ILU(0) preconditioner for $N_p=256,512$. 
For the range $2-512$ of the number of processes, it appears that the computation time 
of the Darcy flow linear system solution remains small  compared with the transport model computation time. 
This can be checked by comparison of the computation times in Table \ref{solver_cpg_time} 
and in Figure \ref{timecpgfrac}. 
This explains the good scalability obtained for both the total and transport computation times 
thanks to the explicit nature of the time integration scheme.

\begin{table}[!htbp]
\caption{Number of linear solver iterations vs. the matrix fracture permeability ratio 
${\Lambda_f\over \Lambda_m}$ for $n_x=128$ and $N_p=2, 128$. }
\label{solver_cpg_perm}
\centering
\begin{tabular}{c|c|c|c|c|c|c}
\hline
& \multicolumn{3}{c|}{$N_p=2$} & \multicolumn{3}{c}{$N_p=128$}\\
\hline
$\Lambda_f/\Lambda_m$ & 20 & 100 & 1000 & 20 & 100 & 1000 \\
\hline
GMRES + Boomer AMG & 15 & 15 & 16 & 15 & 15 & 15  \\
\hline
GMRES + Aggregation AMG & 59 & - & - & 73 & - & -   \\
\hline
GMRES + ILU(0) & 751 & - & - & 633 & - &- \\ 
\hline
\end{tabular}
\begin{tablenotes}
\item -: \small The solver doesn't converge in 1200 iterations.
\end{tablenotes}
\end{table}

\subsection{Fracture network with tetrahedral meshes}

\begin{figure}[!htbp]
\centering
\includegraphics[width=0.38\textwidth]{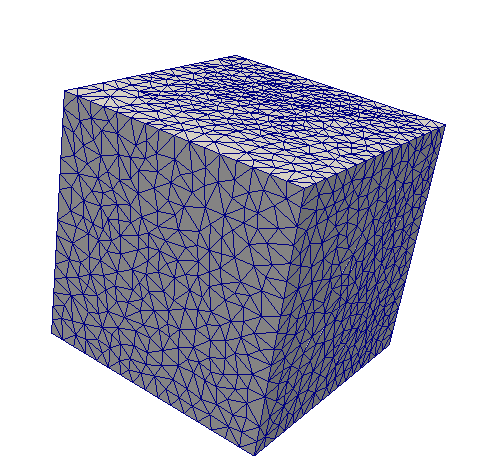}
\hspace{0.05\textwidth}
\includegraphics[width=0.38\textwidth]{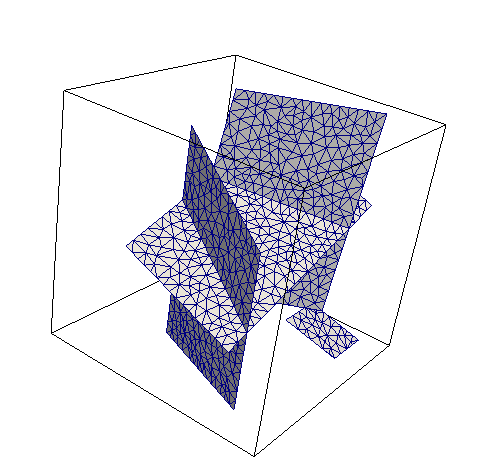}
\caption{Example of tetrahedral mesh of the matrix domain (left) conforming to the fracture network (right).}
\label{meshtet}
\end{figure}

\begin{figure}[!htbp]
          \centering
          \addtocounter{subfigure}{2} 
          \begin{subfigure}{0.48\textwidth}
           \centering
            \includegraphics[width=\textwidth]{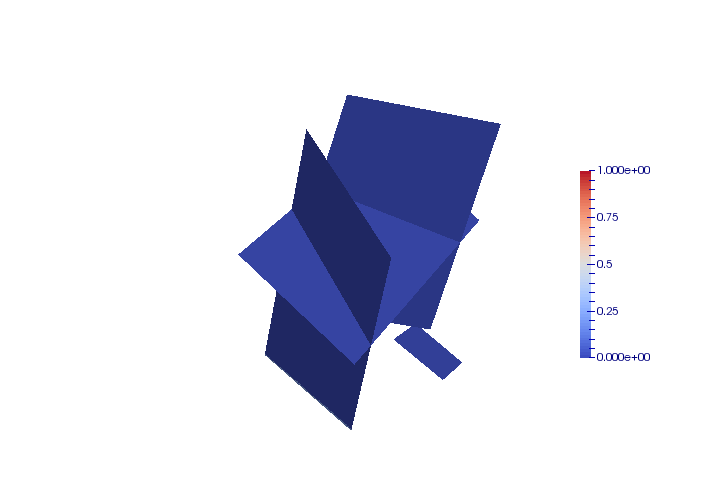}
            \caption{$t=0.0$}
          \end{subfigure}
          \begin{subfigure}{0.48\textwidth}
           \centering
            \includegraphics[width=\textwidth]{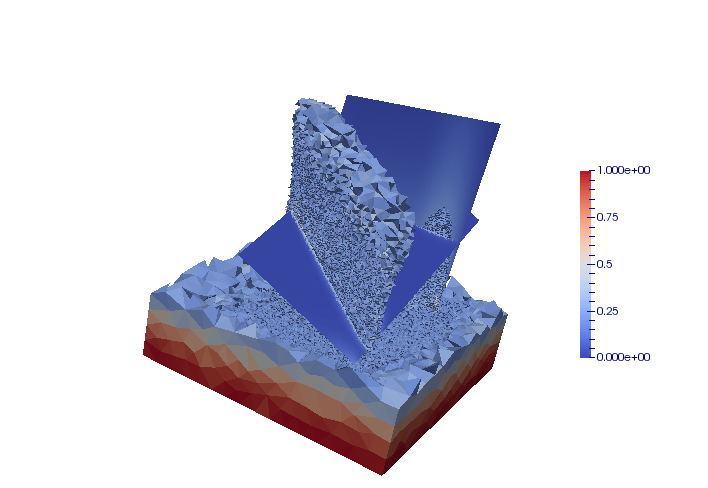}
            \caption{$t=0.2$}
          \end{subfigure}
          \begin{subfigure}{0.48\textwidth}
           \centering
            \includegraphics[width=\textwidth]{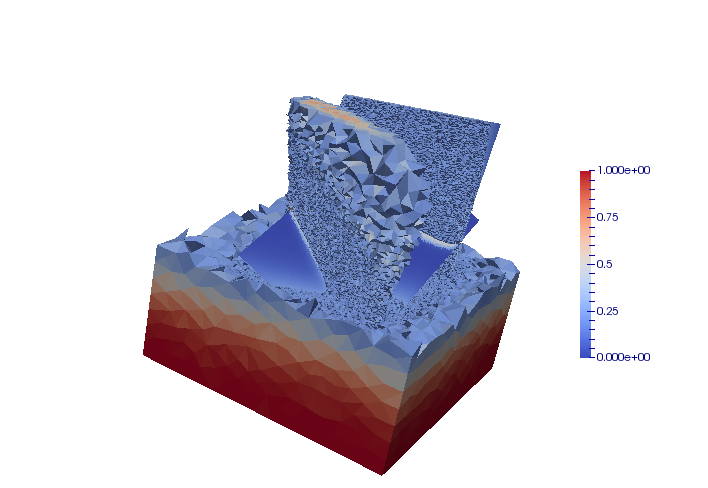}
            \caption{$t=0.4$}
          \end{subfigure}
          \begin{subfigure}{0.48\textwidth}
           \centering
            \includegraphics[width=\textwidth]{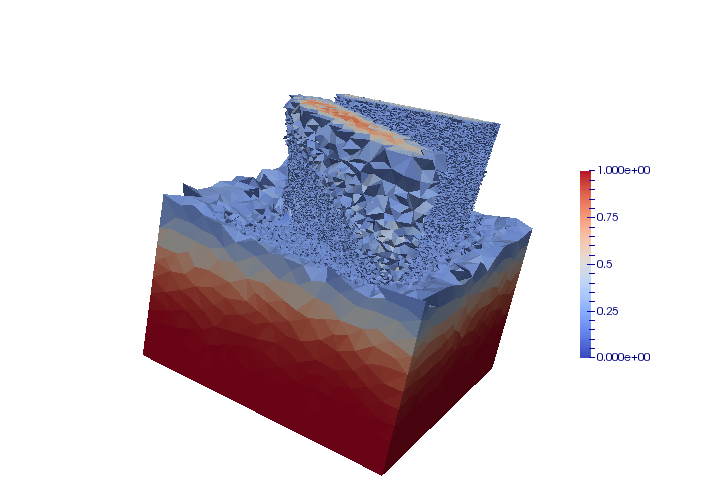}
            \caption{$t=0.5$}
          \end{subfigure}
          \caption{Concentration in the matrix domain and in the fracture network obtained 
                  at different times  for the tetrahedral mesh. 
                   A threshold of 0.2 is used in the matrix domain.}
          \label{simutet}
        \end{figure}

This test case considers tetrahedral meshes of the cubic domain $(0,1)^3$ conforming to the fracture network.  
An example of tetrahedral mesh showing both the matrix domain and the fracture network is exhibited in Figure \ref{meshtet}. All the physical parameters, initial and boundary conditions are the 
same as for the previous test case. 
The mesh used in this subsection contains about $6.2\times 10^6$ cells, $9.7\times 10^5$ nodes and $7.1\times 10^4$ fracture faces. 
Figure \ref{simutet} exhibits the tracer concentrations obtained with this tetrahedral mesh  
at times $t=0$, $t=0.2$, $t=0.4$ and at the final simulation time $t_f = 0.5$.  

As for the previous test case, Tables \ref{solver_tet_nb}, \ref{solver_tet_time} and \ref{solver_tet_perm} 
investigate the performance of the Darcy flow system linear solution for both the GMRES and BiCGStab linear solvers 
and for the same three preconditioners as in the previous test case. The conclusions are basically the same 
as for the hexahedral mesh test case. 
The Boomer AMG preconditioner exhibits an optimal robustness with respect to the matrix fracture permeability ratio  
${\Lambda_f \over \Lambda_m}$. 
On the other hand it requires a rather high number of unknowns per MPI process to 
maintain a good parallel scalability due to the high level of communications in particular in the setup phase. 
The ILU(0) preconditioner can be an interesting alternative but only for 
small size and midly heterogeneous problems. The aggregation AMG preconditioner from the Trilinos library 
used in our test cases 
seems to lack robustness and we did not manage to make it work better through tuning its parameters.

\begin{table}[!htbp]
\caption{
Number of linear solver iterations vs. the number of MPI processes obtained with  different linear solvers and 
preconditioners for the tetrahedral mesh.}
\label{solver_tet_nb}
\centering
\small
\begin{tabular}{c|c|c|c|c|c|c|c|c|c}
\hline
$N_p$ & 2 & 4 & 8 & 16 & 32 & 64 & 128 & 256 & 512 \\
\hline
GMRES + Boomer AMG & 11 & 12  & 12  & 12  & 12  & 12  & 12  & 12  & 12  \\
\hline
GMRES + Aggregation AMG & 38 & 78 & 40 & 39 & 52 & - & 35 & - & 52\\
\hline
GMRES + ILU(0) & 1003 & 725 & 717 & 682 & 667 & 656 & 644 & 629 & 612 \\
\hline
BiCGStab + Boomer AMG & 8 & 7 & 8 & 8 & 8 & 8 & 8 & 8 & 8 \\
\hline
BiCGStab + ILU(0) & 565 & 513 & 527 & 544 & 535 & 483 & 489 & 483 & 473   \\
\hline
\end{tabular}
 \begin{tablenotes}
\small
\item -: The relative residual norm stagnates after a few iterations. 
\item {\color{white} -:}  Some future investigations are necessary.
\end{tablenotes}
\end{table}

\begin{table}[!htbp]
\caption{Linear solver setup phase and solve phase computation times vs. the number of MPI processes 
obtained with  different linear solvers and 
preconditioners for the tetrahedral mesh.}
\label{solver_tet_time}
\centering
\small
\begin{tabular}{c|c|c|c|c|c|c|c|c|c|c}
\hline
$N_p$ & & 2 & 4 & 8 & 16 & 32 & 64 & 128 & 256 & 512 \\
\hline
\multirow{2}{3.1cm}{GMRES \\Boomer AMG}
& Setup time & 12.4 & 7.8 & 4.9 & 5.0 & 4.3 & 6.2 & 7.2 & 13.5 & 22.4 \\
& Solve time & 8.0  & 5.5 & 2.9 & 1.7 & 1.1 & 0.9 & 1.4 & 3.1  & 6.9 \\
\hline
\multirow{2}{3.1cm}{GMRES \\ Aggregation AMG}
& Setup time & 3.7 & 1.9 & 1.2 & 1.8 & 2.1 & 1.6 & 2.9 & 3.3 & 4.7 \\
& Solve time & 19.7 & 20.9 & 5.1 & 2.7 & 2.0 & - & 1.5 & - & 3.0 \\
\hline
\multirow{2}{3.1cm}{GMRES \\ ILU(0)} 
& Setup time & 5.7 & 7.4 & 7.2 & 5.6 & 4.7 & 5.2 & 3.4 & 2.8 & 1.8 \\

& Solve time & 560.6  & 254.4 & 150.0 & 66.5 & 30.1 & 15.2 & 7.7 & 4.1 & 2.8 \\
\hline
\multirow{2}{3.1cm}{BiCGStab \\ Boomer AMG} 
& Setup time & 21.5 & 14.5 & 9.9 & 6.5 & 5.3 & 5.9 & 8.2 & 12.4 & 19.7 \\
& Solve time & 24.0 & 10.2 & 6.1 & 3.5 & 1.8 & 1.5 & 2.1 & 4.3 & 9.5\\
\hline
\multirow{2}{3.1cm}{BiCGStab \\ ILU(0)} 
& Setup time & 5.8 & 6.4 & 6.4 & 5.4 & 4.7 & 5.0 & 3.4 & 2.6 & 1.8  \\
& Solve time & 110.4 & 63.0 & 39.0 & 19.2 & 11.6 & 5.4 & 2.8 & 1.4 & 1.2 \\
\hline
\end{tabular}
\begin{tablenotes}
\small
\item -: The residual norm stagnates after a few iterations. 
\end{tablenotes}
\end{table}

\begin{table}[!htbp]
\caption{
Number of linear solver iterations vs. the matrix fracture permeability ratio 
${\Lambda_f\over \Lambda_m}$ for the tetrahedral mesh and $N_p=2, 128$. }
\label{solver_tet_perm}
\centering
\begin{tabular}{c|c|c|c|c|c|c}
\hline
& \multicolumn{3}{c|}{$N_p=2$} & \multicolumn{3}{c}{$N_p=128$}\\
\hline
$\Lambda_f/\Lambda_m$ & 20 & 100 & 1000 & 20 & 100 & 1000 \\
\hline
GMRES + Boomer AMG &  11 & 13 & 12 & 12 & 13 & 12 \\
\hline
GMRES + Aggregation AMG & 38 & - & - & 35 & - &-  \\
\hline
GMRES + ILU(0) & 1002 & - & - & 644 & - & - \\ 
\hline
\end{tabular}
\begin{tablenotes}
\item -: \small The solver doesn't converge in 1200 iterations.
\end{tablenotes}
\end{table}

\begin{figure}[!htbp]
\centering
\includegraphics[width=0.45\textwidth]{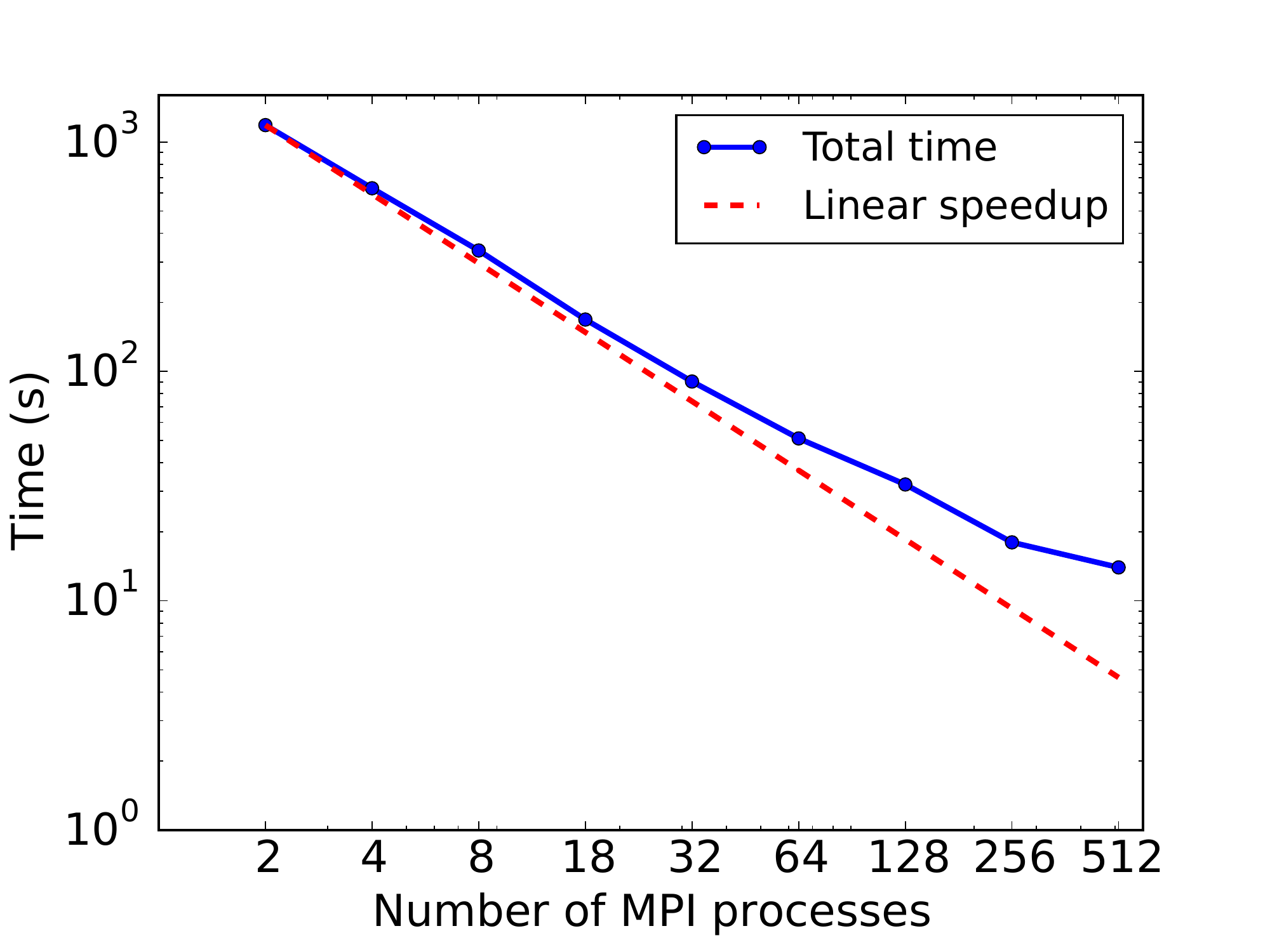}
\hspace{0.05\textwidth}
\includegraphics[width=0.45\textwidth]{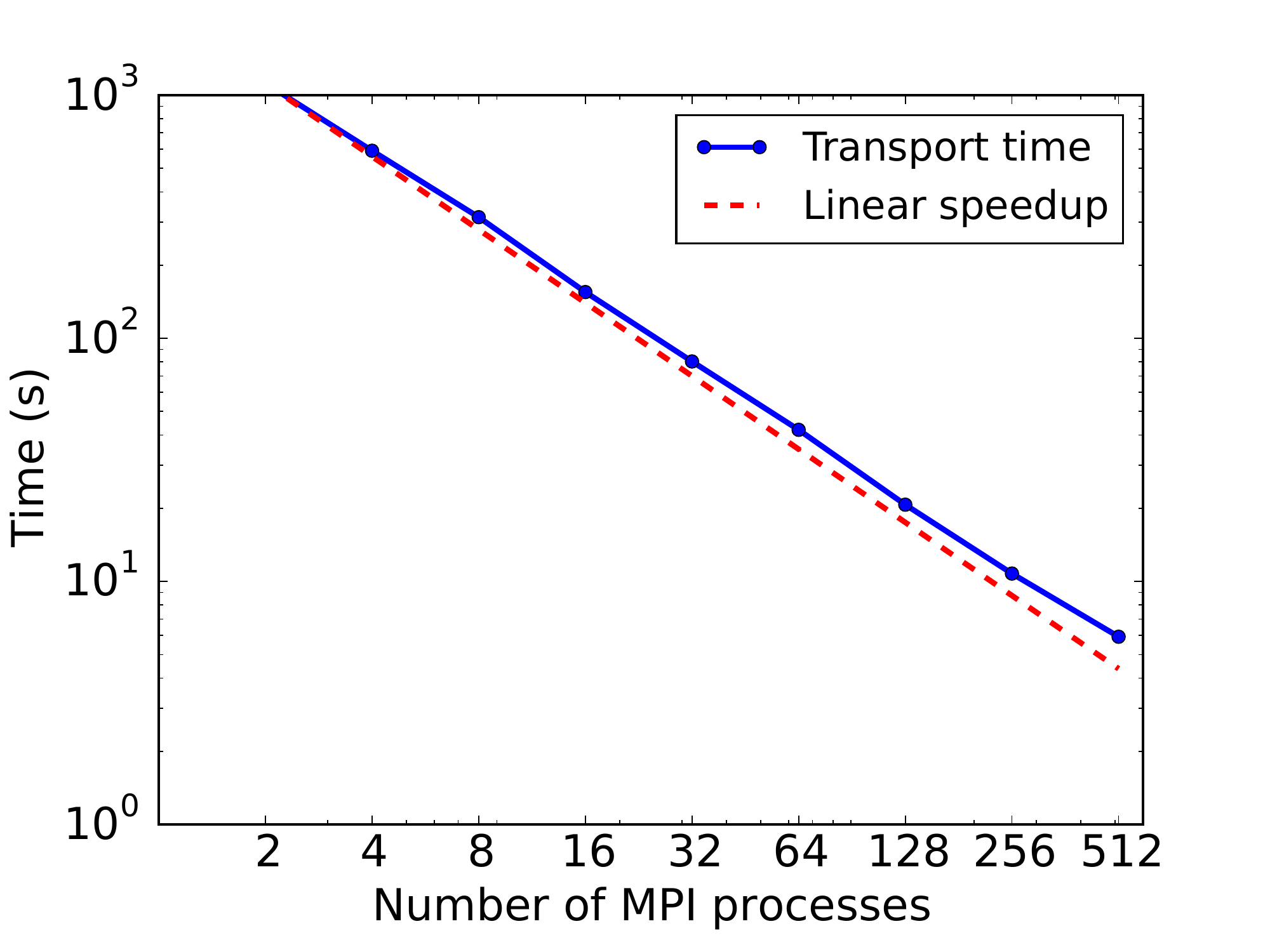}
\caption{Total computation time (left) and computation time for transport model (right)  vs. number of MPI processes with tetrahedral mesh.}
\label{timetet}
\end{figure}

Figure \ref{timetet} plots the total (Darcy flow and transport models) computation time 
and the computation time of the transport model only as a function of the number of processes. 
In these runs the GMRES linear solver is used combined with the 
Boomer AMG preconditioner for $N_p \leqslant 128$ and with the ILU(0) preconditioner for $N_p=256,512$. 
Compared with the previous subsection, an even better parallel scalability 
of the transport model computation time is observed in 
the right Figure \ref{timetet}. This can be explained by the ratio  of roughly $6$ between the number of cells and 
the number of nodes typical of a tetrahedral mesh. For a topologically Cartesian mesh, this ratio is roughly $1$. 
Since the cell concentrations are computed locally in each process, this explains the better scalability observed for 
this tetrahedral mesh compared with the previous hexahedral mesh. 
On the left Figure \ref{timetet}, it is observed that the linear system solution computation time is no longer 
small compared with the transport computation time for $N_p=256$ and $512$.  Hence, it significantly reduces 
the parallel efficiency of the simulation for a large number of processes, say $N_p=256, 512$ in this test case.

\subsection{Application to a complex fracture network}
\label{sec_SimuGGF}

\begin{figure}[!htbp]
\centering
\includegraphics[width=0.5\textwidth]{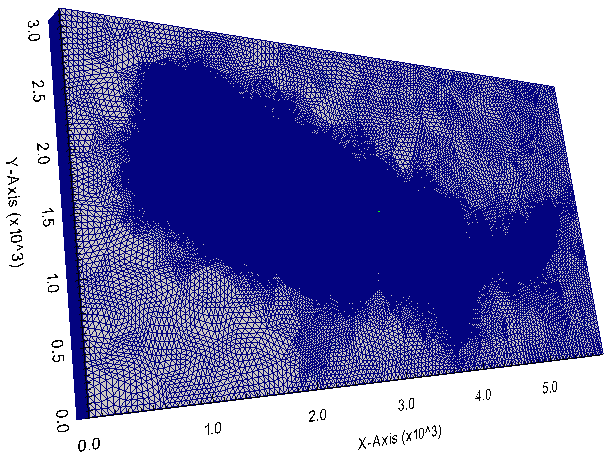}
\caption{Prismatic mesh of the domain $\Omega$ 
defined by the tensor product of a vertical 1D uniform mesh with a 2D triangular mesh.}
\label{mehggfcell}
\end{figure}

\begin{figure}[!htbp]
\centering
\includegraphics[width=0.6\textwidth]{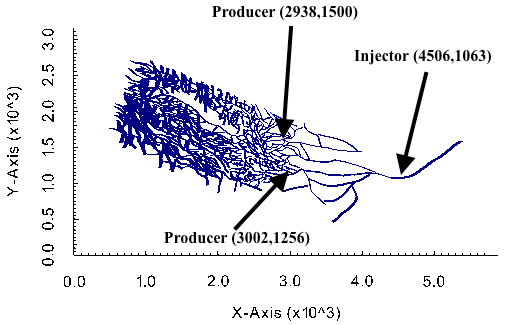}
\caption{Fracture network showing the location of the single 
injection well and of the two production wells.}
\label{meshggffrac}
\end{figure}

In this subsection, our algorithm is applied to a complex fracture network kindly provided 
by M. Karimi-Fard and A. Lap\`ene from Stanford University and TOTAL. 
Figure \ref{mehggfcell} exhibits the mesh of the domain $\Omega=(0,5888.75)\times (0,3157.5) \times (0,250)$ (m) 
which contains about $1.2\times 10^7$ prismatic cells, $6.5\times 10^6$ nodes and $5.13\times 10^5$ fracture faces. 
This 3D mesh is defined by the tensor product of a triangular 2D mesh with a uniform vertical 1D mesh with $24$ intervals.  
The fracture network exhibited in Figure \ref{meshggffrac} contains 581 connected components. It is 
a set of $21376\times 24$ faces of the 3D mesh defined by the tensor product 
of a subset of $21376$ edges of the triangular 2D mesh with the 1D vertical mesh. 
The 2D triangular mesh contains $517540$ cells and is refined in the neighbourhood of the fracture network 
down to an average size of $3.5$ m. 
Figure  \ref{meshggffrac} also shows  the location of the injection well and of the two production wells. 
Each well is vertical of radius $r_w=0.1$ m and its centre 
in the horizontal plane is located at the middle of a fracture edge in the 2D triangular mesh. 
In the vertical direction, only the $12$ fracture faces at the center of the 1D mesh are perforated. 
The permeabilities are isotropic and set to $\Lambda_f=10^{-11}$ ($m^2$) in the fracture network and to 
$\Lambda_m= 10^{-15}$ ($m^2$) in the matrix domain. The porosities are set to $\phi_m = \phi_f = 0.1$, the fracture 
width to $d_f = 1$ m and the fluid viscosity to $\mu = 10^{-3}$ Pa.s$^{-1}$. 

The initial concentration is set to $0$ both in the matrix domain and in the fracture network. 
A total volume of $5.0\times 10^{6}$ m$^3$ is injected in one year at the injector well with a tracer concentration of $1$.  
The pressures of each perforated fracture face $\sigma$ of the producer wells are fixed to $p_w = 0$ and the flow rates  
are given by the Peaceman model 
$$
q_\sigma = WI_\sigma (p_\sigma -p_w),
$$
where $p_\sigma$ is the pressure in the fracture face and $WI_\sigma$ the well index of 
the fracture face. This well index is computed following Peaceman methodology 
\cite{Peaceman78}, \cite{Peaceman83}, \cite{CZ-2009} 
by expanding the fracture face as a box of size $dx\times d_f \times dz$. The analytical pressure solution 
obtained for a vertical well with the well pressure $p_w$, the well radius $r_w$ and the well flow rate $q_w$ per unit length 
is imposed at the 8 corners of the box. Then, the flow rate $q_w dz$ is imposed at the box center 
and the pressure $p_c$ at the box center can be computed analytically using the VAG scheme. We deduce 
the well index $$WI = {q_w dz \over p_c-p_w}$$ leading in this simple case to the analytical formula 
$$
WI = {2 \pi dz \Lambda_f \over  \log(\frac{r_0}{r_w}) }
$$
with 
$$r_0 = D \exp( -{2\pi dz \over C}),$$ and 
$$C = {4 \over 3}( {dx dz \over d_f} + {dx d_f \over dz} + {dz d_f \over dx} ), \quad D = 0.5  \sqrt{ dx^2 + d_f^2 }.$$ 
The production lasts $8$ years.

\begin{figure}
\centering
\includegraphics[width=0.45\textwidth]{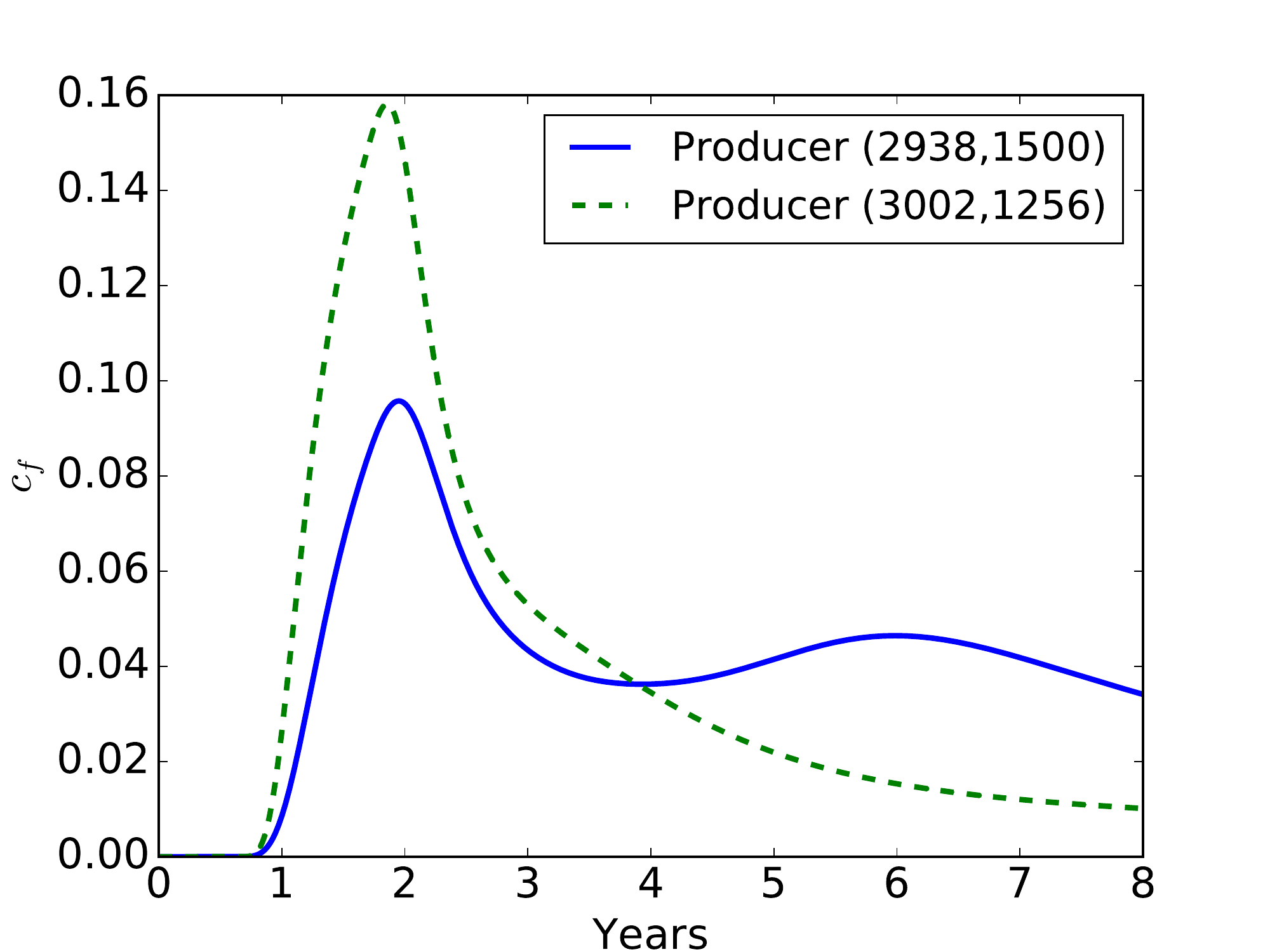}
\hspace{0.05\textwidth}
\includegraphics[width=0.45\textwidth]{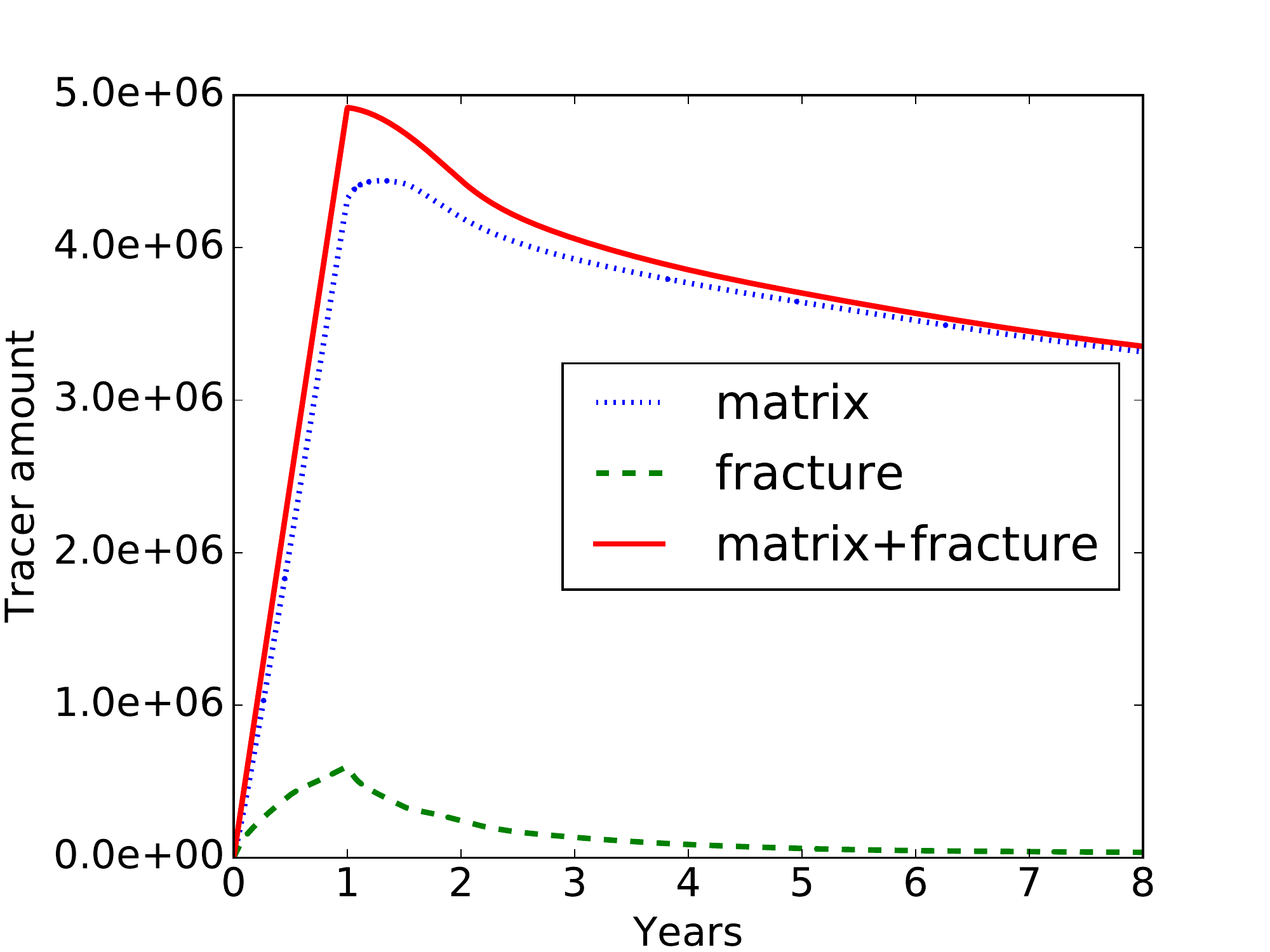}
\caption{Mean tracer concentration in both production wells as a function of time (left)
defined as the ratio between the well tracer flow rate and the well fluid flow rate 
(equal in our case to the well fracture-face tracer concentration). 
Volume of tracer as a function of time in the matrix domain, 
in the fracture network and their sum (right).}
\label{prod10}
\end{figure}

\begin{figure}
\centering
\includegraphics[width=0.5\textwidth]{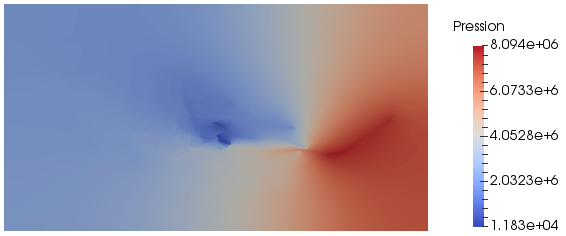}
\caption{Pressure on the matrix domain.}
\label{pressureggf10m}
\end{figure}

\begin{figure}
\centering
\includegraphics[width=0.45\textwidth]{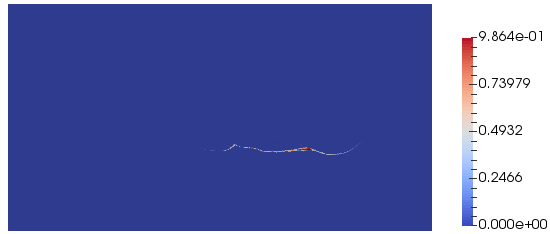}
\hspace{0.05\textwidth}
\includegraphics[width=0.45\textwidth]{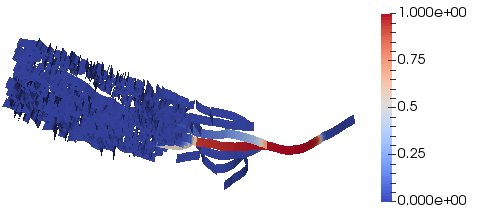}
\caption{Tracer concentration after one year of injection in the matrix domain (left) and in the fracture network (right).}
\label{simuggf_1year}
\end{figure}

\begin{figure}
\centering
\includegraphics[width=0.45\textwidth]{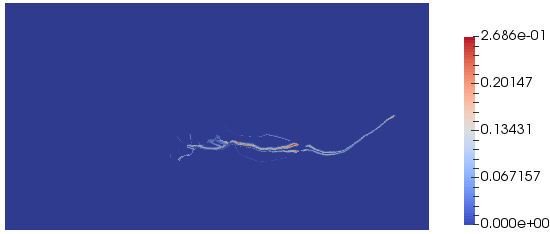}
\hspace{0.05\textwidth}
\includegraphics[width=0.45\textwidth]{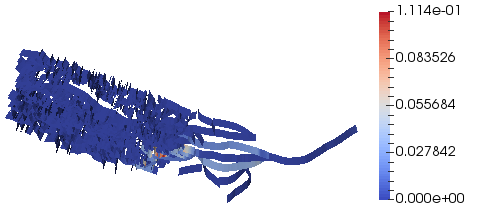}
\caption{Tracer concentration at final time in the matrix domain (left) and in the fracture network (right).}
\label{simuggf}
\end{figure}

\begin{figure}[!htbp]
\centering
\includegraphics[width=0.45\textwidth]{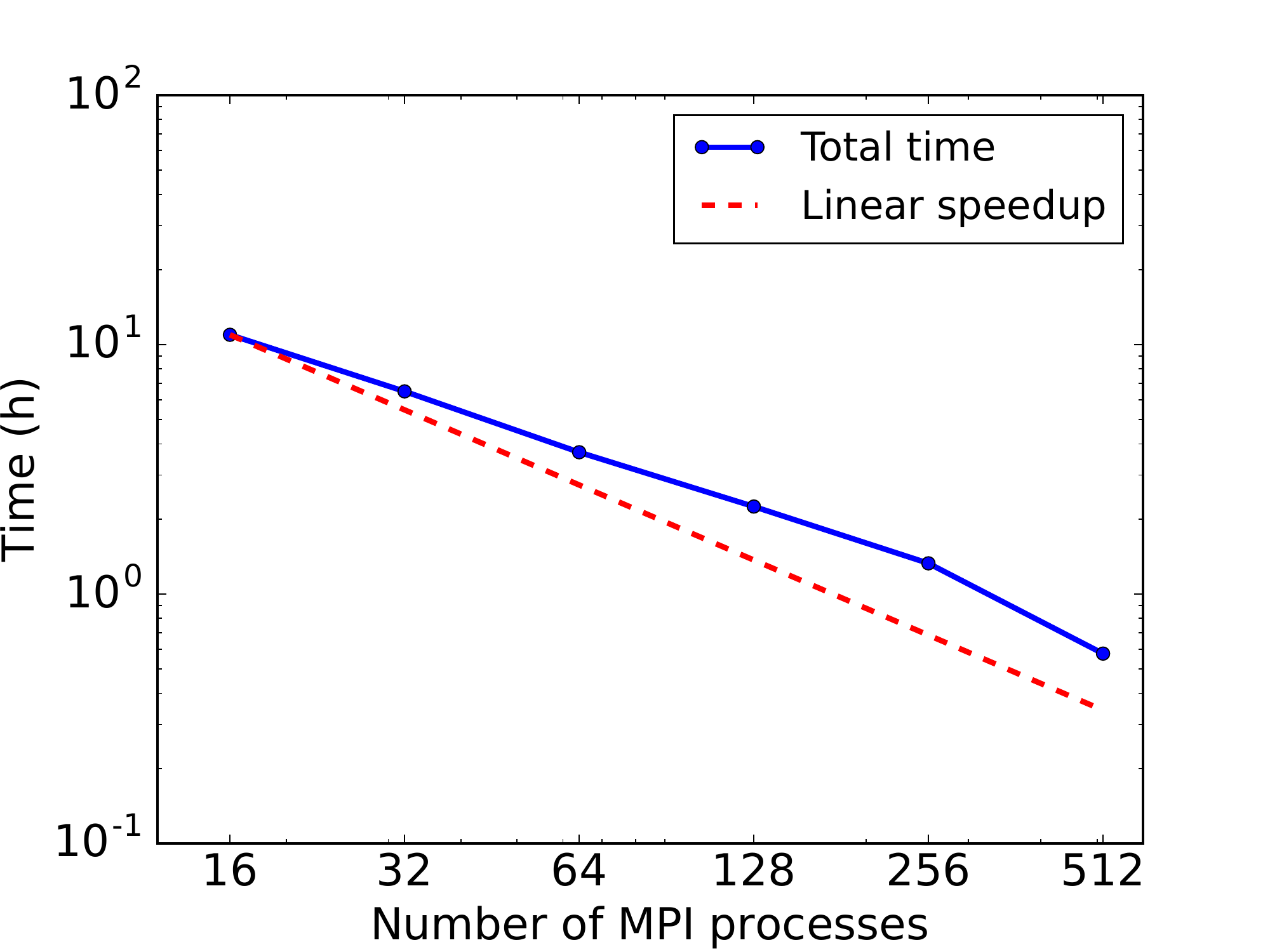}
\hspace{0.05\textwidth}
\includegraphics[width=0.45\textwidth]{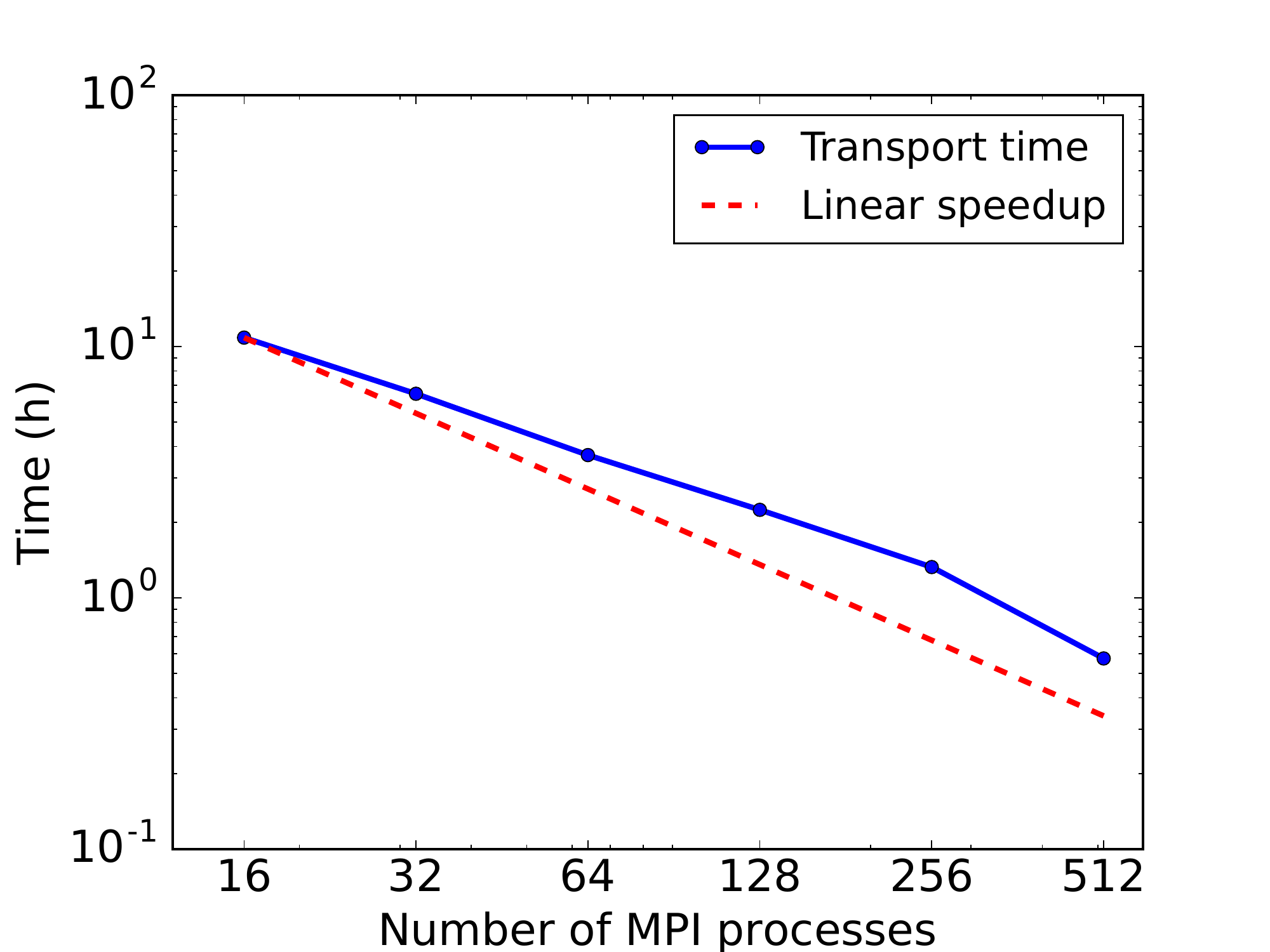}
\caption{Total computation time in hours (left) and computation time for transport model (right)  vs. number of MPI processes with prismatic mesh.}
\label{timeggf}
\end{figure}

Figure \ref{prod10} plots the mean tracer concentration 
in each well as a function of time as well as the total volume of tracer 
as a function of time in the matrix, in the fracture network and their sum. 
Figure \ref{pressureggf10m} exhibit the pressure solution in the matrix domain 
and Figures \ref{simuggf_1year} and  \ref{simuggf} shows 
the tracer concentration after one year of injection and at final time both 
in the matrix domain and in the fracture network. 
Figure \ref{timeggf} shows the total computation times with different number of MPI processes 
$N_p=16,32,64,128,256,512$. It is observed that the total computation time exhibits a rather good scalability. 
In addition, the linear solver (GMRES+Boomer AMG) for the pressure 
converges in no more than 25 iterations whatever the number of MPI processes. Also the comparison of the 
total and transport computation times in Figure \ref{timeggf} shows that the time for the pressure 
solution remains small compared with the transport computation time up to $N_p = 512$. 

\section{Conclusion}
This paper introduced a parallel VAG scheme for the simulation of a hybrid dimensional Darcy flow and transport model  
in a discrete fracture network taking into account the mass exchanges with the matrix.
The convergence of the scheme was validated on two original analytical solutions 
for a flow and transport model that includes fractures. 
The parallel efficiency of the algorithm was studied for different complexities of fracture networks, 
and a large range of matrix fracture permeability ratios and different type of meshes. 
The numerical results exhibit 
a very good parallel strong scalability as expected from the explicit nature of the time integration of the transport model 
with a better result on tetrahedral 
meshes thanks to the communication free computation of the cell unknowns. 
The Darcy flow solution is remarkably robust using the Boomer AMG preconditioner on all types of fracture networks, 
meshes and for all permeability ratios that have been tested.
On the other hand, it requires as usual a rather high number of unknowns per process to maintain a 
good parallel scalability.  
Future work includes the extension of the parallel algorithm to hybrid dimensional multiphase flow models 
and the use of a more accurate second order MUSCL scheme for the transport model.  

\section*{Acknowledgments} 

This work is supported by a joint project between INRIA and BRGM Carnot institutes (ANR, INRIA, BRGM). This work was also granted access to the HPC and visualization resources of ``Centre de Calcul Interactif" hosted by University Nice Sophia Antipolis. This work was partially funded by ADEME BRGM research project Orbou (convention 1305C0131).
We would like also to thank Mohammad Karimi-Fard and Alexandre Lap\`ene for kindly 
providing us the mesh of Section \ref{sec_SimuGGF}.

\end{document}